\documentclass[11pt,preprint]{elsarticle}
\usepackage[toc,title,page]{appendix}
\usepackage{amsmath}
\usepackage{hyperref}
\usepackage{verbatim}
\usepackage{color}
\usepackage{algorithm}
\usepackage{lineno}
\usepackage{multicol}
\usepackage{multirow}
\usepackage[noend]{algpseudocode}
\usepackage{url}
\usepackage{cases}
\usepackage{rotating}
\usepackage{booktabs}
\usepackage{lscape}
\usepackage{subfigure}
\usepackage[acronym]{glossaries}
\usepackage{multirow}
\usepackage{threeparttable}
\newacronym{TW}{tw}{Time Window}
\newacronym{AA}{AA}{Arrival-Arrival}
\newacronym{DD}{DD}{Departure-Departure}
\newacronym{AD}{AD}{Arrival-Departure}
\newacronym{ITSMSM}{ITSMSM}{Integrated Train Scheduling and Maintenance Scheduling Model}
\newacronym{ITSMSOM}{ITSMSOM}{Integrated Train Scheduling and Maintenance Scheduling One-objective Model}
\newacronym{ITSMSODM}{ITSMSODM}{Integrated Train Scheduling and Maintenance Scheduling One-objective Dynamic Model}
\newacronym{STHSRN}{STHSRN}{Space-time high-speed railway network}
\newacronym{TTP}{TTP}{Train Timetabling Problem}
\newacronym{TPP}{TPP}{Train Platforming Problem}
\newacronym{SG}{SG}{Switch Group}
\newacronym{LR}{LR}{Lagrangian Relaxation}
\newacronym{TLSTN}{TLSTN}{Two-level space-time network}
\newacronym{BD}{BD}{Benders Decomposition}
\newacronym{LBBD}{LBBD}{Logic-based Benders Decomposition}
\newacronym{2-L LR}{2-L LR}{Two-level Lagrangian Relaxation}
\usepackage{enumerate}
\usepackage{ulem}
\usepackage{geometry}
\usepackage{natbib}
\bibpunct[, ]{(}{)}{;}{a}{}{,}
\renewcommand\bibfont{\fontsize{10}{12}\selectfont}

\usepackage{amsthm}

\usepackage{comment}
\usepackage{geometry}
\usepackage{caption}

\usepackage{ifthen}
\newcounter{figs}
\newcounter{tabs}
\newcounter{intoappendix}
\usepackage{lineno}

\usepackage{natbib}
\setcitestyle{authoryear}
\biboptions{round,semicolon}
\usepackage{stackrel}
\usepackage{mathtools}
\usepackage{glossaries}
\usepackage{graphicx}
\usepackage{setspace}
\begin{document}

\theoremstyle{plain}
\newtheorem{theorem}{Theorem}[section]
\newtheorem{lemma}[theorem]{Lemma}
\newtheorem{corollary}[theorem]{Corollary}
\newtheorem{proposition}[theorem]{Proposition}

\theoremstyle{definition}
\newtheorem{definition}[theorem]{Definition}
\newtheorem{example}[theorem]{Example}

\theoremstyle{remark}
\newtheorem{remark}{Remark}
\newtheorem{notation}{Notation}
\newcommand*{\qin}{\textcolor{black}}
\newcommand*{\rl}{\textcolor{black}}

\begin{frontmatter}

\title{Solving the train-platforming problem via a two-level Lagrangian Relaxation approach}

\author[First]{Qin Zhang}
\ead{17114203@bjtu.edu.cn}
\author[Second]{Richard Martin Lusby}
\ead{rmlu@dtu.dk}
\author[First,Third]{Pan Shang\corref{cor1}}
\ead{shangpan@bjtu.edu.cn}
\author[Fourth]{Chang Liu}
\ead{18811443345@163.com}
\author[First]{Wenqian Liu}
\ead{19114044@bjtu.edu.cn}

\address[First]{School of Traffic and Transportation, Beijing Jiaotong University, Beijing 100044, China}
\address[Second]{Department of Technology, Management and Economics, Technical University of Denmark, Kgs. Lyngby 2800, Denmark}
\address[Third]{Key Laboratory of Transport Industry of Comprehensive Transportation Theory (Beijing Jiaotong University), Ministry of Transport, Beijing 100044, China}
\address[Fourth]{Transportation \& Economics Research Institute, China Academy of Railway Sciences Corporation Limited, Beijing
100081, China}

\cortext[cor1]{Corresponding author}

\begin{abstract}
High-speed railway stations are crucial junctions in high-speed railway networks. Compared to operations on the tracks between stations, trains have more routing possibilities within stations. As a result, track allocation at a station is relatively complicated. In this study, we aim to solve the train platforming problem for a busy high-speed railway station by considering comprehensive track resources and interlocking configurations. A two-level space-time network is constructed to capture infrastructure information at various levels of detail from both macroscopic and microscopic perspectives. Additionally, we propose a nonlinear programming model \rl{that minimizes a} weighted sum of total \qin{travel} time and total deviation time for trains at the station. We apply a Two-level Lagrangian Relaxation (2-L LR) to a linearized version of the model and demonstrate how this induces a decomposable train-specific path choice problem at the macroscopic level that is guided by Lagrange multipliers associated with microscopic resource capacity violation. As case studies, the proposed model and solution approach are applied to a small virtual railway station and a high-speed railway hub station located on the busiest high-speed railway line in China. Through a comparison of other approaches that include Logic-based Benders Decomposition (LBBD), we highlight the superiority of the proposed method; on realistic instances, the 2-L LR method finds solution that are, on average, approximately 2\% from optimality. Finally, we test algorithm performance at the operational level and obtain near-optimal solutions, with optimality gaps of approximately 1\%, in a very short time.
\end{abstract}


\end{frontmatter}
\section{Introduction}

High-speed railway transportation systems have gained significant popularity in recent years due to being fast, safe, punctual, and environmentally friendly. As a result, high-speed railway networks worldwide have rapidly expanded, with substantial increases in both the total length of railway lines and the number of trains. The comparatively limited growth in the number of high-speed railway stations, however, has resulted in the need for each station to accommodate a greater number of railway lines. For instance, the Nanjingnan railway station, one of the main hub stations in the Chinese high-speed railway network, currently connects five railway lines, and it is expected that even more will need to be accommodated in the future. As a result, train routing at this type of station can become very complicated. Therefore, it is important for train planners and dispatchers to enhance the utilization 
of existing station track resources.

The train platforming problem, which is also known as the track allocation problem, is typically solved with the input of a train timetable that specifies train arrival and departure times at the station. This problem optimizes the allocation of track resources at the station to each train in the timetable while satisfying station regulations and safety requirements. Variations of this problem from the strategic perspective focus on designing station layouts, see e.g., \citet{zwaneveld1970decision,qi2016integrated}. At the tactical level, the track allocation to trains based on either fixed arrival and departure times in a predefined timetable or adjustable arrival and departure times in a flexible timetable is solved~\citep{lusby2011railway}. The fixed variant can also be viewed as a flexible problem with zero flexibility. Such problems are offline problems. However, at the operational level, the problem focuses on re-assigning conflict-free routes to trains in a short time-frame, making it an online problem. To the best of our knowledge, the train platforming problem is often solved manually in practice~\citep{zhang2020overview,dang2021low}. In this paper, we consider the train platforming problem at the tactical level for a busy hub station that connects two or more railway lines. The arrival and departure times of trains are allowed to change within a certain range. Modifying the arrival and departure times of trains at one station can lead to changes in the train schedules at all stations in the interconnected railway network where the station is located. The feedback from an individual station can be used to adjust platform plans for other stations in the railway network. However, in this study, we focus on the track allocation at a single busy station within the railway network. Optimizing the train platforming at busy hub stations is seen as being more important than simpler intermediate stations located on a single railway line due to the intricate nature of the station infrastructure that must be considered and the higher volume of train traffic (from multiple different lines) that must be coordinated. \rl{The feasibility of the train platform plan for a large station therefore takes priority over those of small stations for connected railway lines
~\citep{burggraeve2017robust}.} Additionally, we also attempt to solve the problem at the operational level by generating a feasible track allocation plan within a very limited solving time when disruptions occur, such as train delays or track failure.

Unlike the simple parallel railway lines between two stations, the layout of a high-speed railway hub station is complicated and has a number of critical locations such as crossings or switches. Therefore, the precision with which station infrastructure is modelled is important when addressing the train platforming problem since this determines how accurately the track utilization can be described~\citep{pellegrini2014optimal}. A train route at a station connects a sequence of several critical track resources. The complex station layout enables trains to choose multiple reasonable routes in the station and dwell on or pass through the corresponding platform tracks. This has resulted in a significant increase in possible conflicts among train routes at a station. To ensure safety, the occupation and release times of station resources for each train route are controlled by a complex interlocking system. All track resources that comprise a route are simultaneously locked and then incrementally released in the route-locking, sectional-release configuration. Such a system has been applied in many countries, such as China, German and Japan. To optimize train operations at a station, an appropriate representation of infrastructure should be developed to comprehensively illustrate the occupation and release times of track resources by trains and intelligently reduce inefficient paths in such an interlocking system. Typically, the train platforming problem refers to scenarios in which a single route connects a station's entrance or exit point to a platform designated for train operations. The scenario involving numerous route choices between station's entrance and exit points considering critical switches or small track unit is referred to as the train routing problem~\citep{lusby2011railway}. 

Although the track structure within a station may be intricate, in practice, only one basic route that connects one entering (leaving) point and a specific platform is assigned to a certain train. This is because many impractical routes are not considered in an interlocking system, even if they are physically feasible~\citep{zhang2020simultaneously}. Therefore, we focus on the train platforming problem. We assume that only one route exists to/from a platform from/to any station boundary point.
However, the proposed approach can be modified to address the train routing problem since we consider the detailed utilization of switches.

Modeling approaches based on space-time networks have recently been applied to optimize the operations of high-speed railway systems. We refer readers to~\citet{zhang2020simultaneously,zhang2022heuristic}, \citet{meng2014simultaneous} and \citet{zhang2020joint} for additional information. In these works, a single-level space-time network is typically constructed for a specific railway line. Such a single-level network can easily represent the occupation of track resources by trains in the time dimension. However, the size of the network is extremely sensitive to the precision with which time and infrastructure are modelled, and this
 poses major challenges for modeling and solution algorithms. To overcome this, multi-level frameworks have been proposed to address train operation optimization problems, as demonstrated in studies such as \citet{goverde2016three} and \citet{wang2023joint}. However, these approaches often involve coupling full models at different levels with different focuses, which maintains the large size of the lower-level model. Therefore, we suggest to circumvent these limitations by using a novel two-level space-time network that facilitates a detailed modelling of 
 infrastructure (i.e., microscopic level) and couples this with a less detailed representation (macroscopic level) when solving the train platforming problem.  This paper therefore makes the following contributions to the literature on train-platforming:

\begin{enumerate}
    \item We develop a two-level space-time network that captures important infrastructural information with different precision levels. The topology of a station is modeled at a macroscopic route level and detailed track resources (i.e., switches in the station) are modeled at a microscopic level. This representation facilitates the consideration of the route-locking, sectional-release interlocking configuration in the train platforming problem, and reduces overall network size. This network is suitable for addressing problems where a decision is made to address the issue from an aggregated perspective, while simultaneously considering disaggregated resources.
    \item Based on the two-level network, we formulate a 
    two-level model that minimizes a weighted sum of the total \qin{travel} time and total \qin{shift} time from the desired time for all trains. Unlike existing multiple-level models, in which certain constraints from otherwise independent full models are often coupled, our proposed model operates as a cohesive framework, with each level responsible for distinct constraints. Complex capacity constraints are only enforced at the microscopic level with detailed precision for 
    infrastructure, whereas other operational constraints are guaranteed at the macroscopic level with a coarser granularity. Additionally, the linking constraints used to connect the decision variables associated with the respective levels facilitate the representation of interlocking configurations.
    \item We devise a \gls{2-L LR} algorithm to solve the linearized mathematical model. 
    This allows us to decompose the problem into a master problem that comprises a set of macroscopic train-specific blocks, which can each be solved using a dynamic programming approach, and a subproblem that assesses microscopic feasibility. Lagrange multipliers that are associated with the microscopic level are first aggregated and then used in the macroscopic level to provide information on track resource overutilization. 
    \item  We compare the performance of the proposed algorithm to that of \rl{a} \gls{LBBD} method and a lazy constraint method on a small virtual station, as well as a Chinese high-speed railway station over a time horizon of a full day with a 15-second time discretization and a total of 287 trains. The results demonstrate that the proposed algorithm can provide high-quality solutions to the train platforming problem at the tactical level. Additionally, we can obtain a high-quality feasible solution within a very short time-frame when solving the online problem. 
\end{enumerate}

The remaining sections of this paper are organized as follows. Section~\ref{sec:literature} provides a brief review of the most relevant literature. In Section~\ref{sec:description}, we present the physical track resources 
at both the macroscopic and microscopic levels, and describe the construction of the two-level space-time network. The proposed nonlinear programming formulation is presented in Section~\ref{sec:model}. Section~\ref{sec:alg} linearizes the two-level model and introduces the proposed \gls{2-L LR} algorithm. Results from numerical experiments conducted on a small virtual station and a Chinese high-speed railway station are reported in Section~\ref{sec:experiment}. Finally, Section~\ref{sec:conslusion} concludes the paper and provides insights for future research directions.

\section{Literature review}
\label{sec:literature}

The train platforming problem is a classical problem in the railway planning and rescheduling process. This problem is an NP-complete problem~\citep{cardillo1998k}, and various models and algorithms have been developed over the past few decades to solve it. \citet{lusby2011railway} reviewed work on common approaches to train platforming and routing problems from strategic, tactical, and operational perspectives, including conflict graph methodology, constraint programming, heuristics, and alternative graph formulation. \citet{cardillo1998k} first defined the train platforming problem as a graph coloring problem and proposed a heuristic method to solve it. \citet{rodriguez2007constraint} presented a constraint programming model that made use of a simulator to address the track allocation problem at railway junctions 
without platforms. The author describes a branch-and-bound strategy to effectively solve this problem. \rl{To avoid the computational burden of solving a mathematical model, 
\citet{carey2003scheduling} developed a heuristic approach that mirrors manual scheduling procedures to generate a feasible platform plan for a large complex railway station.  
This is then extended by~\citet{carey2007scheduling} to consider a network of stations.}
\citet{billionnet2003using} formulated two integer programming models based on \citet{cardillo1998k} and solved them directly using an integer programming solver. In subsequent research, linear programming formulations with various methodologies have been widely employed (e.g., \citet{chakroborty2008optimum,caprara2011solution,lusby2011routing}). To solve the online train routing at complex junctions, \citet{pellegrini2014optimal} developed a mixed integer programming model that minimized the maximum and total secondary delay suffered by trains. A rolling-horizon solution framework was proposed to obtain results. Based on this, \citet{sama2016ant} further designed an ant colony algorithm to address the real-time train routing selection problem given all train routing possibilities at the station. \citet{d2007branch} first introduced the concept of alternative graphs, proposed by \citet{mascis2002job}, to optimization problems associated with railway operations. This formulation can effectively describe operational time constraints and the incompatible occupation of railway resources by trains and has become a popular modeling approach for the train 
\qin{platforming} problem \citep{corman2009rescheduling}. 

Among the various models proposed for solving this problem, the most challenging and crucial issue is the representation of the conflict-free occupation constraint for station resources. 
\citet{zwaneveld1996routing} and \citet{zwaneveld2001routing} preprocessed conflicts between pairs of train space-time routes and formulated the routing problem as a node packing problem that included train-pair-route-dependent conflict clique constraints. \citet{caprara2011solution} investigated a general formulation with a quadratic objective function and linearized it by introducing auxiliary variables and constraints. The resulting integer linear programming model enforces safety requirements by forbidding the simultaneous assignment of incompatible patterns based on clique inequalities. A similar formulation was also adopted in \citet{lamorgese2016optimal}, where trains were assumed to require the same running time to traverse a station. However, the authors acknowledged that this assumption was not appropriate for addressing train platform problems at large-scale stations. \citet{cacchiani2015tutorial} defined a pattern incompatibility graph to enumerate incompatible train-pattern pairs implicitly. Additionally, corresponding clique constraints were constructed to prevent assigning pairwise incompatible patterns. Similarly, \citet{lu2022train} proposed a mixed integer programming model with the big-{\textit M} method to enforce train-pair-route-dependent conflict-free constraints for the train platforming and rescheduling problem. The occupation times of routes with different interlocking configurations were preprocessed by an ordered conflicting route pair calculation method at the aggregated route level. \qin{The authors proposed a heuristic method that cannot alone provide a lower bound to evaluate the solution quality. The optimality gap for the heuristic \rl{that considers a} route-locking, sectional release configuration is \rl{obtained by comparing the quality of the solutions obtained with the best bounds found by applying CPLEX to a mathematical model that considers a sectional-locking, sectional-release configuration.}} As the problem size increases 
such a modelling approach may lead to a large number of conflicting train-route pairs and corresponding capacity constraints. Furthermore, slight changes in the possible routes of trains can result in significant changes in constraints, as conflicts need to be recomputed. To overcome the problem, \citet{caimi2011new} proposed the use of a microscopic resource tree to represent all potential train paths within a station. They formulated the train routing problem as a resource-constrained multicommodity flow model, incorporating strong cutting planes based on the maximal potential conflict clique. However, it was noted that the explicit preprocessing of train paths in the microscopic level encountered memory issues prior to problem-solving. This motivates the idea of storing information at different levels of detail. \citet{wu2013track} used \rl{a} big-{\textit M} formulation to ensure that trains occupying the same station tracks satisfied the minimum safety time interval. However, the minimum time requirements for using the same switches were not considered. \rl{Motivated by this, \citet{kang2019stochastic} considered the problem of providing a more balanced track allocation for large railway stations in a stochastic setting. Tracks in the interlocking area were modelled as several groups of switches, and, from a safety perspective, constraints that ensured their conflict-free allocation were included.} 
\citet{sels2014train} analyzed the resource occupation times of train movements in a station with a rough route-locking and route-release interlocking system. A mixed integer programming model with Boolean variables was constructed based on \citet{billionnet2003using} to maximize the number of trains assigned to a suitable platform. \citet{bevsinovic2019stable} assumed that all train orders were fixed and proposed a two-level multi-objective optimization model to route trains through a railway station based on an event-activity network. To guarantee the feasibility of a platform plan, the authors enforced the minimum headway by directly preprocessing the blocking times 
in a max-plus automata model at the microscopic level and then simulated the results at the macroscopic level. \rl{Robustness of complex railway stations was also considered by \citet{burggraeve2017robust} and included in the form of buffer time between two trains that required the same microscopic track resource.} \rl{The authors proposed a non-linear model that was subsequently linearized and solved with CPLEX.} \citet{lusby2013set} formulated a set packing model and developed a branch-and-price algorithm based on a discretized space-time microscopic resource constraint system to address train re-routing problem under disruption at a complicated junction. The train flow and infrastructure constraints were enforced only on a microscopic level. Note that the utilization of track resources at a railway station can be directly affected by the interlocking system. Therefore, \citet{corman2009rescheduling} presented a quantitative comparison of two interlocking systems for an aggregated formulation (route-locking, route-release) and disaggregated formulation with non-sequence-dependent setup times (sectional-locking, sectional-release). Results showed that these formulations provide upper and lower bounds for the disaggregated formulation with sequence-dependent setup times (route-locking and sectional-release). The authors emphasized the need for an effective formulation and efficient algorithm to analyze the potential of the interlocking system quantitatively with a configuration of route-locking and sectional-release, which is commonly used in practice. 

The train timetabling problem focuses on optimizing the arrival and departure times of trains at each station while ensuring compliance with safety requirements and the efficient use of resources. Therefore, some researchers have incorporated 
routing decisions at stations for individual trains into their train timetables. \citet{carey1994model} addressed a train timetabling problem that incorporated platform choice for a general railway network, without considering the interlocking systems of stations. \citet{petering2016mixed} formulated an optimization model for a unidirectional railway line with intermediate stations to solve an integrated cyclic timetabling and platforming problem. The big-{\textit M} formulation was used to enforce the headway restriction on tracks between stations and safety constraints at each station. \citet{pellegrini2019efficient} constructed a mixed integer linear programming model to minimize train delay by optimizing train local routing and timetabling problems based on a microscopic infrastructure representation. A branch-and-bound solution procedure with valid inequalities was proposed to solve this model and several realistic instances verified the algorithm's efficiency. A microscopic model was proposed in \citet{zhang2019microscopic} to optimize train timing, sequencing and routing decisions, and the timing of maintenance tasks in an integrated manner. The interlocking system considered was similar to route-locked and route-release. Near-optimal solutions were obtained by a heuristic algorithm that repeatedly solved train scheduling, routing, and maintenance tasks. \citet{zhan2015real} and \citet{zhan2016rolling} investigated an optimization model from a macroscopic perspective to solve the real-time periodic train timetabling problem for a high-speed railway network. The headway requirements at stations were implicitly constrained by the maximum number of trains that a station could accommodate at a certain time instant. Platform assignment for a specific train and the occupation of detailed track resources within a station were not considered. \citet{bevsinovic2016integrated} proposed a macro-micro framework to obtain a conflict-free, robust, and stable timetable. 
The authors transformed the results of feasibility analysis in the microscopic model by dynamically tightening or relaxing the headway and running time constraints in the macroscopic model. The headway constraints were enforced by preventing multiple prepared paths in a single conflict path clique from being scheduled simultaneously. \citet{van2017designing} solved the train timetable adjustment problem to find a feasible alternative timetable by retiming, reordering, canceling, and short-turning trains when maintenance was required at a station or on the tracks between stations. The authors defined the number of available tracks at a station at a certain time using a large set of auxiliary variables to guarantee that at least one station track would be available for incoming trains. The station details, particularly the interlocking system, however, were not adequately formulated within this integrated problem.

Space-time networks have been adopted to optimize railway operations. \citet{ambrosino2021innovative} explained how to model the capacity requirements for problems involving time and space decisions in a general operation-time-space network and presented an application for scheduling port rail shunting operations. \citet{zhang2020simultaneously} proposed a binary integer programming model and two decomposition algorithms based on a mesoscopic space-time network to optimize train timetabling and platforming problems simultaneously for a high-speed railway line.  In a subsequent work, \citet{zhang2022heuristic} further expanded their investigation by incorporating maintenance schedules and designed a dynamic time window algorithm to generate train timetables, track allocation plans, and derive maintenance schedules for a high-speed railway network. To generate a timetable considering both passenger and railway infrastructure perspectives, \citet{xu2021train} developed a constrained multi-commodity network flow formulation to integrate train timetabling and platforming plans while considering passenger behavior on a single, one-way railway line. A Lagrangian relaxation-based heuristic was developed to solve this problem and instances constructed based on a practical high-speed railway line verified the effectiveness of the proposed method. When the time horizon is finely discretized, the size of a space-time network will \qin{increase dramatically}~\citep{zhang2022heuristic}. 
Therefore, strategies \rl{must} be defined to reduce network and related model sizes~\citep{ambrosino2021innovative}. \qin{Recently, \citet{zhang2023optimal} simultaneously optimized the scheduling of locomotive operations and the platforming of railcars at conventional passenger railway stations.} \rl{The proposed model was based on a time-space-state network in which the capacity requirements were explicitly enforced using the state dimension.} 
\qin{However, one disadvantage of this approach is that creating multiple states for each space-time arc results in a large-scale network. Furthermore, using a time-discretization of one-minute
for the occupation of a track unit is not accurate enough.} 

Multiple-level approaches that aim to formulate railway transportation optimization problems have attracted more and more attention from researchers. In \citet{zhu2014scheduled}, a three-level space-time network was developed to tackle the service network design problem for freight railway transportation. The three levels describe the consolidation flow with the cars, ranging from single-car units to blocks and trains. The capacity constraints were enforced at each level. \citet{goverde2016three} proposed a three-level framework for railway timetabling problems to integrate timetable construction and performance evaluation problems in an iterative process. Three different objectives were optimized by the three corresponding independent models. The microscopic level within this framework was designed to accurately calculate the running time and blocking time for train trajectories. Similarly, \citet{wang2023joint} implemented a coupled multi-resolution space-time railway network to jointly optimize the train timetabling and routing problems for a single railway line with simple railway stations. The macroscopic level optimized the train schedule by considering two boundary nodes and one aggregated platform node for each station. The microscopic level considered the flow balance requirements and assigned appropriate routes for the given macroscopic path. After the decomposition of the full model, the macroscopic model was solved with additional constraints that were iteratively generated from the microscopic model. \citet{xu2020integrated} formulated an integrated micro-macro nonlinear mixed-integer programming model to incorporate train speed control into the railway timetabling problem. Due to its complexity, the optimization model was reformulated based on a single-level space-time-speed network in which segment-level train trajectories were presolved. After relaxing the capacity constraints, a final solution was generated by detecting and resolving conflicts in the resulting solution. The authors focused on train scheduling between stations and  overlooked the utilization of station resources by trains. A similar transformation is applied in \citet{zhang2020simultaneously} and \citet{zhang2022heuristic}, where the occupations of station track resources are implicitly integrated into the top-level arcs by preprocessing incompatible arcs when jointly optimizing the train timetabling and platforming problems. To address the complexity of the problems while maintaining the feasibility of platform plans, the authors considered a route-locking and route-release interlocking system at the station level, where all track sections comprising a station route were assumed to be released simultaneously. However, this precision may not be suitable for specialized train platforming problems, potentially leading to inefficient utilization of station capacity. 
To better demonstrate the contribution of this paper, Table ~\ref{literature} provides a summary of related research. For each paper, we list the level of the constructed network, the specific problem perspective, the modeling method used, the solution algorithm, the interlocking configuration applied and the largest problem. For the largest problem, we list the time precision in minutes (m), seconds(s), the number of trains, and the time horizon in hours (h), minutes (m), days (d).

\section{Problem description}
\label{sec:description}
\begin{figure}[h!]
    \centering
    \includegraphics[scale=1.05]{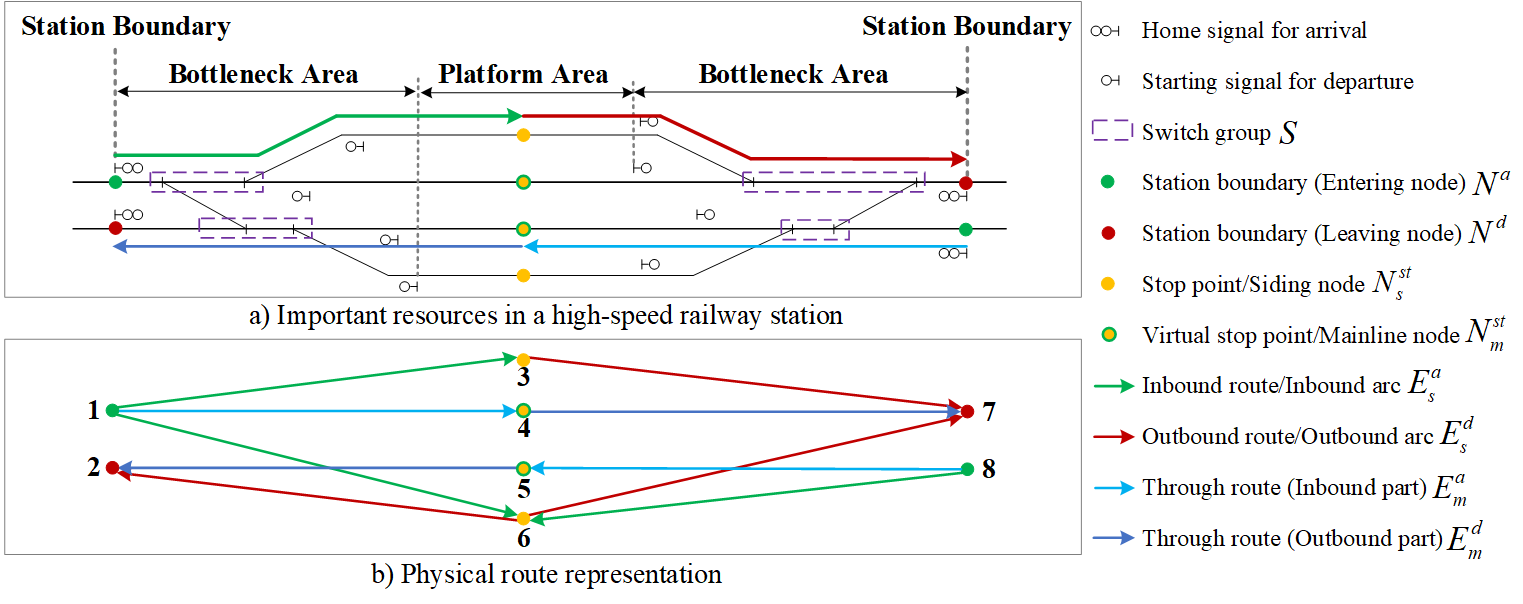}
    \caption{\rl{An overview of important 
    track resources within a station and the representation of the physical route.}}
    \label{fig:main_resource}
\end{figure}

Track resources at a railway station consist of tracks in the platform area and tracks in the {\it bottleneck areas} that lead to and from the platform area, as shown in Fig.~\ref{fig:main_resource}-a). Tracks in the platform area comprise {\it siding tracks}, for the dwelling of stopping trains, and {\it mainline tracks}, for the passing of nonstop trains. Trains occupy track resources in the bottleneck area using different types of train routes~\citep{zhang2022heuristic}. Stopping trains travel from the station boundary along {\it inbound routes}, dwell at the stop points on siding tracks for passengers to board and alight, and then leave the station along {\it outbound routes} after the desired dwelling time. Nonstop trains traverse the mainline tracks via {\it through routes}. As in
~\citet{zhang2022heuristic}, we introduce virtual stop points on the mainline tracks and separate the routes of nonstop trains into inbound and outbound components. Compared to 
the platform area, resource allocation in the bottleneck area is more complex. An interlocking system controls the switches and signals on a specific train route and prevents conflicting movements between trains. \qin{Usually, certain switches cannot be \rl{simultaneously assigned to two trains due to infrastructure requirements.} In some cases, it may be physically possible to simultaneously assign some switches to different trains, but this is not allowed for safety reasons.}   
\rl{Such switches can be termed a {\it \gls{SG}}}. \qin{\rl{Furthermore}, if the running time between two switches that belong to different track sections is very short, then these switches can also be referred to as a single \gls{SG}.} \rl{In this paper, we consider \glspl{SG} as important resources in the bottleneck area rather than individual switches. This helps  to reduce the problem size; however, \qin{the level of detail of the representation is defined by the modeler.} \rl{If,} necessary, it is still possible to address a switch-level problem by treating each individual switch as its own \gls{SG}.}

Most existing research on the train timetabling problem defines the granularity of infrastructure from two main perspectives: macroscopic and microscopic. From a macroscopic perspective, stations are simplified either as points or as two separate points for arrivals and departures, and the railway track segments that connect two stations are treated as unified entities. On the other hand, the microscopic perspective focuses on a more deta-

\renewcommand{\arraystretch}{1.5}
\begin{landscape}
\begin{table}[htbp]
\centering
\caption{Summary for the related publications.}
\label{literature}
\small{
\begin{tabular}{p{100pt}p{22pt}p{58pt}p{197pt}p{110pt}p{67pt}p{71pt}}
\toprule
Publication &Level& Perspective & Modelling method &Algorithm&Interlocking&Largest problem\\
\midrule
\citet{zwaneveld1996routing}&Single&Routing&Preprocessed train-pair-route-deviation-based; Node packing &Branch-and-cut-based&RLSR&1m/27/1h\\
\citet{zwaneveld2001routing}&Single&Routing&Preprocessed train-pair-route-based; Weighted node packing &Branch-and-cut-based&RLSR&-/79/1h\\
\citet{caimi2011new}&Single&Routing&Preprocessed resource-based; Conflict graph&Branch-and-cut-based&RLSR&-/67/1h\\
\citet{lusby2013set}&Single&Routing&Preprocessed  resource-based; Node packing&Branch-and-price-based&RLSR&15s/66/2h\\
\citet{zhang2019microscopic}&Single&Routing&Resource-based &Heuristic&RLRR&1s/58/-\\
\citet{sama2016ant}&Single&Routing&Preprocessed route-based; Assignment model& Ant colony& RLSR& -/40/1h\\


\citet{bevsinovic2019stable}&Two&Routing&Preprocessed resource-simulate-based model&Greedy heuristic&-& 1s/14/30m\\
\citet{wang2023joint}&Two&Routing&Resource-based; Coupling two full models&LR\&CG&SLSR&1s/48/5h\\
\citet{goverde2016three}&Three&Routing&Three separated modules with feedback&Iterative&-&5s/41/1h\\
\citet{pellegrini2014optimal}&-&Routing&Resource-based&CPLEX\& Heuristic&RLSR\&RLRR&1s/45/1h\\
\qin{\citet{burggraeve2017robust}}&\qin{Single}&\qin{Routing}&\qin{Resource-based}&\qin{CPLEX}&\qin{RLRR}&\qin{-/85/1h}\\
\qin{\citet{kang2019stochastic}}&\qin{Single}&\qin{Platforming}&\qin{Resource-based}&\qin{Simulated annealing}&\qin{-}&\qin{1m/28/4h}\\
\citet{xu2021train}&Single&Platforming&Route-based&LR&-&1m/150/1d\\%
\citet{zhang2020simultaneously}&Single&Platforming&Route-based&LR; ADMM; CPLEX&RLRR&1m/40/380m\\
\citet{zhang2022heuristic}&Single&Platforming&Route-based&Heuristic&RLRR&30s/23/325m\\
\citet{lu2022train}&-&Platforming&Route-based&Heuristic&RLSR\&RLRR \&SLSR&\qin{1s/100/4.5h}\\
\qin{\citet{zhang2023optimal}}&\qin{Single}&\qin{Platforming}&\qin{Resource-based }&\qin{LR;ADMM}&\qin{RLSR}&\qin{1m/40/200m}\\
\textbf{This paper}&\textbf{Two}&\textbf{Platforming}&\textbf{Route-resource-based single model}&\textbf{Two-level LR; LBBD}&\textbf{RLSR\&RLRR}&\textbf{15s/287/1d}\\
\bottomrule
\multicolumn{7}{l}{LR: Lagrangian relaxation; CG: Constraint generation;
ADMM: Alternating direction method of multipliers; RLRR: Route-locking, route release; } \\
\multicolumn{7}{l}{RLSR: Route-locking, sectional-release; 
 SLSR: Sectional-locking, sectional-release; -: Not applicable }
\end{tabular}}
\end{table}
\end{landscape}
\renewcommand{\arraystretch}{1}

\noindent iled representation of track resources, such as switches within a station or block sections between stations. Our previous studies (\citet{zhang2020simultaneously,zhang2022heuristic}) focused on the simultaneous optimization of train timetables and platform assignments. In these studies, we considered each platform track within a station as an individual point, while still treating the tracks between stations as a whole. This representation effectively captured station routes and was referred to as the mesoscopic level. However, this paper specifically focuses on the station resources themselves. A station cannot therefore be reduced to a single point from a modelling perspective. Within the scope of this paper, we use the term macroscopic to refer to the coarser representation of the station infrastructure and note that this corresponds to the mesoscopic level that was defined in our previous work. 

\subsection{Physical infrastructure representation}

To represent train movements at 
a station, we describe all track resources using a two-level network with different infrastructure 
details from macroscopic and microscopic perspectives.

At the macroscopic level, the physical layout of the station is represented by a directed graph $G^P=(N,E)$, similar to~\citet{zhang2022heuristic}. The node set $N$ comprises four different subsets of node types: the entering node set $N^a$, the leaving node set $N^d$, the siding node set $N^{st}_{s}$, and the mainline node set $N^{st}_{m}$. \qin{The mainline nodes can only \rl{be used by} nonstop trains \rl{as} no trains can dwell on the mainline tracks \rl{of a} station.} The physical arc set $E$ includes inbound arc set $E^a$ and outbound arc set $E^d$, which respectively, represent the inbound and outbound routes. \hspace{0.5em}The inbound arc set $E^a$ can be further classified into arc sets $E^a_s$ and $E^a_m$,\hspace{0.5em} to distinguish between arcs that connect a station boundary with a stop point from those that connect a station boundary with a virtual stop point. $E^d=E^d_s\cup E^d_m$ is defined similarly. For each physical arc $e=(i,j)\in E$, we define the origin node as $o_e=i$, the destination node as $d_e=j$, and the running time on arc $e$ as $t_e$. Fig.~\ref{fig:main_resource}-b) illustrates how we model station infrastructure at the macroscopic route level.


At the microscopic level, we consider \glspl{SG} as essential track resources of routes in the bottleneck area. We use $S$ to denote the  set of all \glspl{SG}. The track resources in the platform area at the microscopic level are modelled using the macroscopic level node sets $N_s^{st}$ and $N_s^{m}$.

There is a connection between the occupation of track resources in the bottleneck area in the macroscopic and microscopic representations. For each physical inbound or outbound arc $e\in  E^a \cup E^d$ in the macroscopic level, we define all \glspl{SG} that are associated with $e$ as $\phi(e)$.

\subsection{Two-level space-time network representation}

Optimizing the train platforming at a busy high-speed railway station that enforces a route-locking and sectional-release interlocking system requires careful consideration of the infrastructure granularity when using space-time networks to model train movements. At the microscopic level, train path conflicts can be detected by considering each specific track resource~\citep{caimi2011new}. However, relying solely on a microscopic 
\qin{representation} results in an excessive number of space-time arcs and poses challenges, especially at stations with complex track layout, in describing the route-locking configuration. Additionally, this representation may result in an excessively large solution space since it could allow for the inclusion of impractical routes, which would need to be thoughtfully refined. On the other hand, describing the station infrastructure solely from the macroscopic level introduces challenges in formulating the sectional-release configuration.  This is because the utilization of switch groups in the bottleneck area becomes train-pair-route dependent, adding complexity to the methodologies.  However, at the macroscopic level, the node and arc sets are relatively small, making it easier to enforce the basic operational requirements. We therefore propose a two-level space-time network to capture important infrastructural details from both macroscopic and microscopic perspectives:

\begin{itemize}
    \item The macroscopic level is based on the station route level and is used to describe the macroscopic trajectory of the trains, describing their entry at specific station boundaries, their passage along designated platform tracks, and their subsequent exit via distinct station boundaries. This level enforces basic operational constraints, such as flow balance and any dwelling requirements, and reflects the sequence of critical track resources associated with a specific train route. 
    However, the capacity requirement is neglected. Moreover, this level can effectively capture the synchronized locking time of the switch groups of one route.
    \item The microscopic level is used to describe the occupation time of microscopic track resources and to enforce safety requirements. At this level, the sectional-release configuration can be readily executed for individual switch groups that correspond to a macroscopic route. As a result, enforcing the safety requirements is straightforward since we can guarantee the exclusive occupation of all microscopic resources.
\end{itemize}

\subsubsection{Macroscopic space-time network}


Based on the macroscopic physical network $G^P$, we can construct a macroscopic space-time network as a weighted directed graph $G^T=(V,A)$ to represent train movement along track resources at a railway station. However, the detailed occupation of microscopic resources is overlooked. Time is discretized for the entire planning horizon and is given by the set of time periods $T=\left\{1,2,...,|T|\right\}$. The vertex set $V$ includes a virtual source vertex $\overline{o}$, a virtual sink vertex $\overline{k}$
\qin{, the entering vertex set $V^a$, the leaving vertex set $V^d$, the siding vertex set $ V^{st}_s$ and the mainline vertex set $V^{st}_m$. Each vertex $(i,t)\in V^a\cup V^d \cup V^{st}_s\cup V^{st}_m$ is obtained by copying the corresponding physical node $i\in N^a\cup N^d \cup N^{st}_s\cup N^{st}_m$ at each time \rl{period} $t\in T$.}
The arc set $A$ consists of the source arc set $A^o$, the sink arc set $A^k$, the arrival arc set $A^a$, the departure arc set $A^d$ and the siding waiting arc set $A^{st}$. \qin{A source arc $(i,j,t,\tau)
\in A^o$ connects the virtual source vertex $(i,t)=\overline{o}$ to an entering vertex $(j,\tau)\in V^a$, and a sink arc $(i,j,t,\tau)
\in A^k$ connects a leaving vertex $(i,t)\in V^d$ to the virtual sink vertex $(j,\tau)=\overline{k}$.} These two types of arcs represent the start and end of the train flow. \qin{An arrival arc $(i,j,t,\tau)\in A^a$ starts from an entering vertex $(i,t)\in V^a $ and ends at a siding vertex $(j,\tau)\in V^{st}_s$ or mainline vertex $(j,\tau)\in V^{st}_m$. Such an arc represents train movement \rl{along an} inbound route $(i,j)\in E^a$. A departure arc $(i,j,t,\tau)\in A^d$ connects a siding vertex $(i,t)\in V^{st}_s$ or mainline vertex $(i,t)\in V^{st}_m$ to a leaving vertex $(j,\tau)\in V^d$. Such an arc represents train movement \rl{along an} outbound route $(i,j)\in E^d$. A siding waiting arc $(i,i,t,t+1) \in A^{st}$ connects two neighboring vertices $(i,t)\in V^{st}_s$ and $(i,t+1)\in V^{st}_s$ corresponding to the same siding node $i\in N^{st}_s$ and represents the dwelling time on a siding track.} To guarantee the feasibility of a solution, trains that cannot be assigned to actual platform tracks can be assigned a virtual path. The virtual path \qin{$g_p$} starts from the source vertex \qin{$\overline{o}$} and ends at the sink vertex \qin{$\overline{k}$}. The weight of \qin{the} virtual path represents the cost of cancellation. \qin{For simplicity, we use $v$ with $n_v=i$, $t_v=t$ to denote the vertex $(i,t)\in V$ and $g$ with $o_g=i, d_g=j, \hat{o}_g=t, \hat{d}_g=\tau$, $e_g=(i,j)$ to represent arc $(i,j,t,\tau)\in A$.}  \rl{We denote the set of trains as $F$. Similarly, we denote the set of available vertices for train $f\in F$ as $V_f\subseteq{V}$ and the set of available space-time arcs for train $f\in F$ as $A_f=A^o_f\cup A^k_f\cup A^a_f\cup A^d_f\cup A^{st}_f \cup \left\{g_p\right\}\subseteq{A}$.}

\begin{figure}[h!]
    \centering
    \includegraphics[scale=0.8]{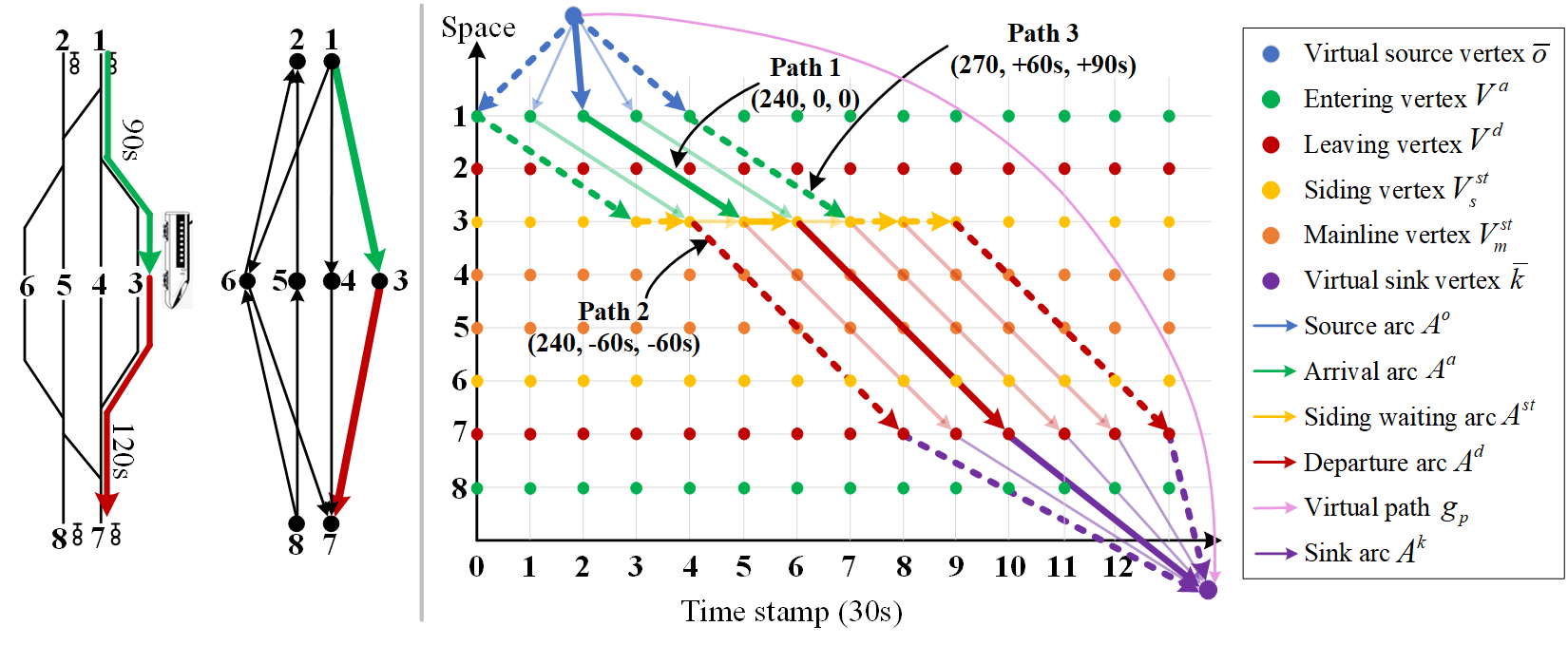}
    \caption{\qin{Macroscopic space-time network representation.}}
    \label{fig:meso_space-time net}
\end{figure}

\qin{Fig.~\ref{fig:meso_space-time net} \rl{provides an overview} of the macroscopic space-time network. To \rl{aid the} visualization, only a \rl{subset of the arcs for one train are shown.}} The train is allowed to dwell at Node 3, and the desired arrival and departure times for this train are assumed to be 150 and 180 seconds, respectively. Meanwhile, the maximum arrival and departure shifts are set to 60 and 90 seconds. We show three possible space-time paths of 
\qin{the} train. 
\qin{Two of these are dashed and one is in bold.} For each path, we list the \qin{travel} time, shift in arrival time, and shift in departure time. One can see that Path 1 is a desirable path because its shifts of arrival and departure times are both zero. The \qin{travel} time of Path 2 is equal to that of Path 1, but the arrival and departure times of Path 2 are 60 seconds earlier than the desirable times. The \qin{travel} time of Path 3 is greater than that of the other two paths. Specifically, Path 3 arrives at the station 60 seconds later than expected with an additional 30 seconds of dwelling and eventually departs from the station 90 seconds later than the desired time. If none of these paths conflict with the paths of other trains, then Path 1, is clearly the optimal choice.

\subsubsection{\qin{Space-time representation for the microscopic level}}

\qin{The space-time representation for the resources at the microscopic level is designed to represent the occupation of specific resources, i.e., \gls{SG} resources in the bottleneck area and platform resources along the time dimension by trains.} At this level, the set of space-time track resources is defined as $B= B^s \cup B^{st}$, where $B^s$ and $B^{st}$ represent the sets of space-time \gls{SG} and siding track resources, respectively. For \qin{a physical} \gls{SG} resource $s\in S$, the corresponding microscopic space-time resource set is denoted as $l=(\qin{s, }s, t, t+1)\in B^s$ for $t=1,2,…,|T|-1$. For a siding 
\qin{track resource} $i\in N^{st}_s$, the microscopic space-time resource set is defined as $l=(i, i, t, t+1)\in B^{st}$ for $t=1,2,…,|T|-1$. 
\qin{We use $s_l$ to denote the corresponding \gls{SG} for the microscopic \rl{space-time} resource $l\in B^s$ and $n_l$ to denote the corresponding siding track resource for the microscopic \rl{space-time} resource $l\in B^{st}$. The time period \rl{associated with resource $l\in B$ is denoted as $\hat{o}_l=t$. The set of microscopic space-time resources for train $f$} is defined as $B_f=B_f^{s}\cup B_f^{st}$.} \qin{It is worth mentioning that the time granularity for the resources at the microscopic space-time level \rl{need not be the same as that} at the macroscopic space-time level,}  
\rl{and we discuss how to address such a case in the presentation of our approach.} 
The siding waiting arcs \qin{at the macroscopic level} only represent the actual occupations of siding tracks by trains, whereas we consider additional factors \qin{at the microscopic level}, which will be described in the following subsections.



\subsubsection{\qin{Connection between the macroscopic and microscopic levels}}

\vspace{5pt}
 \noindent {\bf (1) Track resources in the bottleneck area}
\vspace{5pt}

Before a train physically occupies an inbound or outbound route, some time should be allocated for route formulation while taking into account the safety distance~\citep{pellegrini2014optimal}. We refer to this as the {\it formulation time}. All the \qin{switch} resources that belong to a route are reserved within the formulation time prior to being traversed. 
The corresponding \qin{switch}  resources of one route can be utilized by other trains a certain time after the \qin{switch} has been traversed. The duration of this additional time is termed the {\it clearing time}. \qin{In the literature, two main infrastructure utilization modes in the bottleneck area have been applied and modeled depending on the configuration of the interlocking system at a high-speed railway station.} 

 $\bullet$ {\bf Mode 1. Route-locking, route-release.}
 
All \qin{switch} resources associated with one route are locked and released simultaneously and can be reserved by another train after the actual occupancy of the entire route and clearing time. 
 
 $\bullet$ {\bf Mode 2. Route-locking, sectional-release.}
 
All \qin{switch} resources associated with one route are locked simultaneously. However, they are released incrementally and can be reserved by another train after the actual occupancy of the resources and clearing time. 

\begin{figure}[h!]
    \centering
    \includegraphics[scale=0.69]{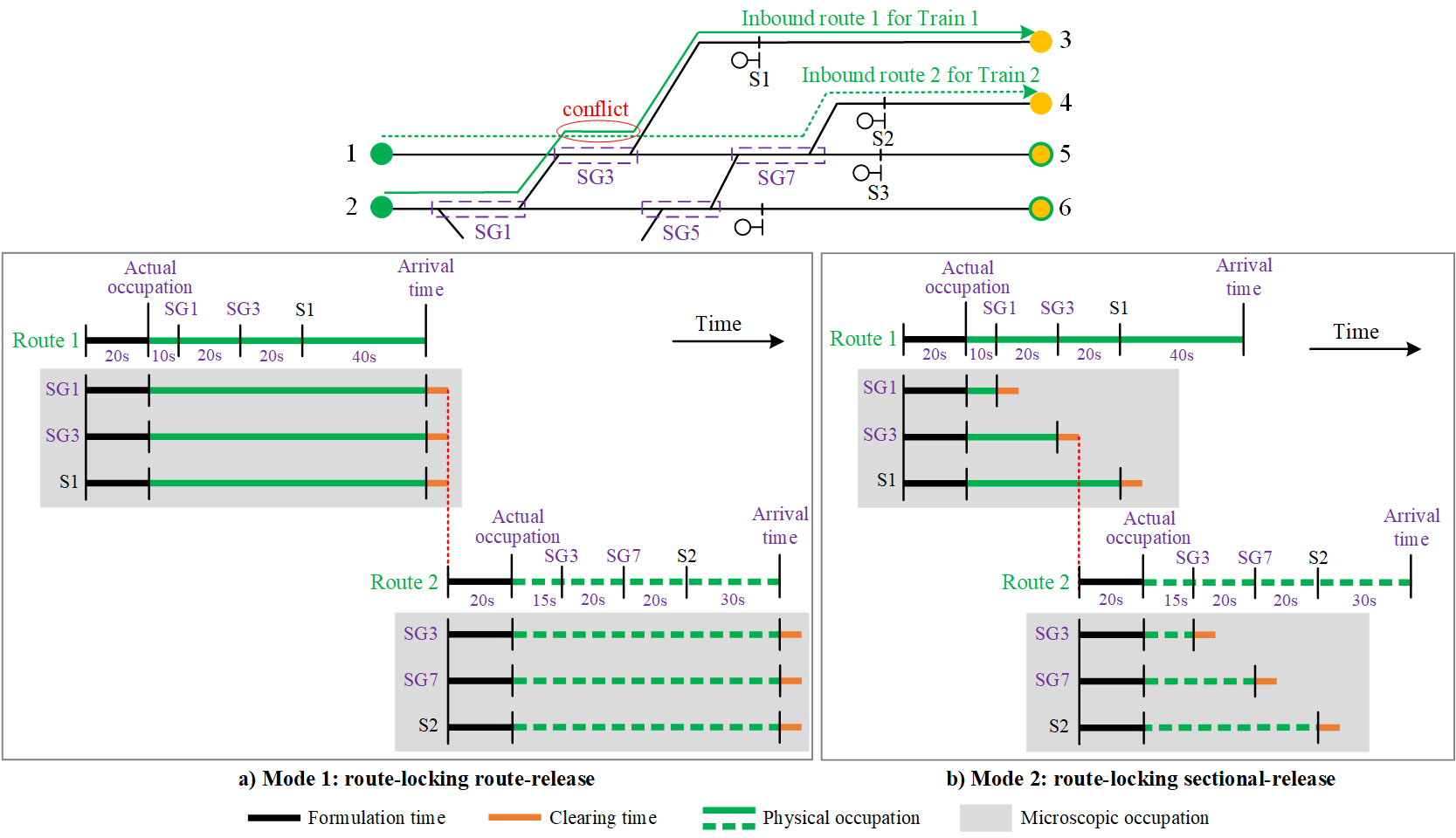}
    \caption{Occupation of resources at two levels for the arrival movement.}
    \label{fig:map_arrival}
\end{figure}

\qin{In practice, the route-locking, sectional-release mode is more widely adopted. Route-release can be considered as sectional-release \rl{when the release times of all switch resources of the route is the same}. Therefore, the configuration of route-locking and sectional-release is considered in the remainder of this paper. As we described earlier, we considered \gls{SG} resources rather than individual switch in the bottleneck area.} Fig.~\ref{fig:map_arrival} presents a simple example to highlight the differences between these two modes \rl{when considering} \glspl{SG}. Train 2 can occupy the inbound route up to 60 seconds earlier in the route-locking, sectional-release mode compared to the route-locking, route-release mode. A similar observation can be made when trains are leaving a station.  For the sake of simplicity, the minimum headway required for two trains to occupy the same inbound or outbound route and corresponding microscopic resources is defined as the sum of the formulation time and clearing time.

\vspace{5pt}
 \noindent {\bf (2) Track resources in the platform area}
\vspace{5pt}

\begin{figure}[h!]
    \centering
    \includegraphics[scale=0.96]{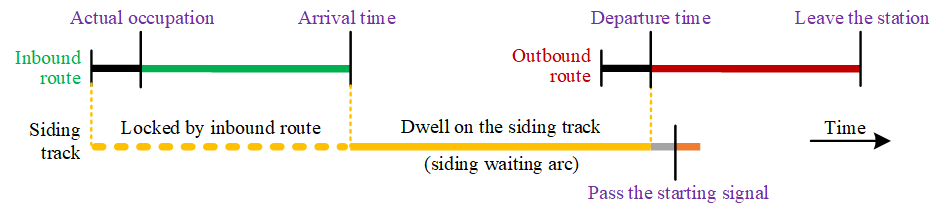}
    \caption{Occupation of siding tracks.}
    \label{fig:occupation_siding}
\end{figure}

The occupation times of the siding tracks are illustrated in Fig.~\ref{fig:occupation_siding}. A train reserves the siding track the moment its corresponding inbound route is prepared. After dwelling on the siding track for a predetermined time, the train departs at its departure time. In practice, this train finishes its occupation of the siding track when it passes the starting signal. The siding track can then be assigned to other trains after the clearing time has elapsed. For the sake of simplicity, the 
headway, which is the minimum time required for two consecutive trains to occupy the siding track, is defined as the time interval between the departure time of the first train and the siding track's release time. In other words, the end time of one actual occupation of a siding track is assumed to be the departure time of the train. 

\vspace{5pt}
 \noindent {\bf (3) Connection between track resource occupation at two different levels}
\vspace{5pt}

\begin{figure}[h!]
    \centering
    \includegraphics[scale=1.04]{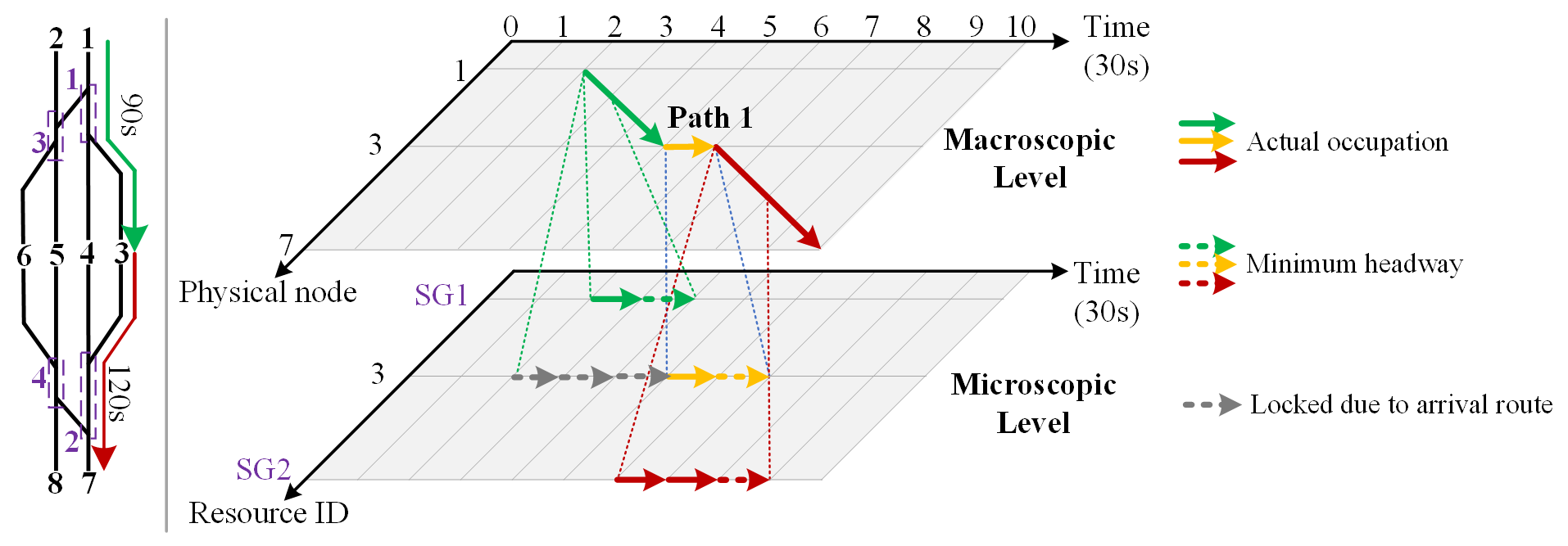}
    \caption{\qin{Connection between resource occupation at two different levels.}}
    \label{fig:map_two_ST_net}
\end{figure}

We describe the actual occupation of infrastructure by trains at the macroscopic level and represent the occupation of resources, including the formulation time, clearing time, and locked time, with finer precision in the microscopic level. Fig.~\ref{fig:map_two_ST_net} illustrates the connection between the occupations of track resources by Path 1 in Fig.~\ref{fig:meso_space-time net} for the proposed 
\qin{two levels. For each arrival and departure arc $g\in A^a\cup A^d$, we define its associated microscopic \gls{SG} resources as $\phi_g^{SG}=\left\{l|l\in B^{s}, s_l\in \phi(e_g),\hat{o}_l\in[\hat{o}_g,\hat{o}_g+T_{e_g}^{s_l}+\Delta r)\right\}$, where $T_{e}^{s}$ represents the occupation time of the {SG} resources $s\in S$ for the physical route $e$ in the given interlocking configuration and $\Delta r$ denotes the minimum headway required for the usage of the same \gls{SG} resource. The set of microscopic siding resources that will be simultaneously reserved for the arrival arc $g\in A^a$ is defined as $\phi_{g}^{st}=\left\{l|l\in B^{st},o_{l}=d_g,\hat{o}_{l}\in [\hat{o}_g,\hat{d}_g)\right\}$. 
\rl{The occupation time of} 
microscopic resources in these two sets is determined by the number of time \rl{periods}. Therefore, both \rl{sets must} be carefully constructed \rl{given the} relationship between the actual time and the time granularity \rl{respectively used} in the macroscopic and microscopic levels, especially when the time granularity in \rl{each level} is different.}

\section{Mathematical formulation}
\label{sec:model}

We present a \gls{TLSTN} model to solve the train platforming problem in a busy high-speed railway station given the station's interlocking configuration.

\subsection{Assumptions}
\label{subsec:assump}

The following assumptions are made when modeling this problem:
\begin{itemize}

    \item All infrastructure information is available. This includes the station layout, the running times on the inbound and outbound routes, and the minimum headway between two train occupations on the same track resource.
    \item All train data is known. This includes the desired arrival and departure times at the (virtual) stop points as well as the maximum allowable \qin{shift} from these times for all trains in the ideal timetable. The ideal timetable in this paper refers to the timetable that does not take into account the track allocation in stations at the tactical level or a feasible timetable before disruption happens at the operational level.
    \item \qin{The running time on \rl{any}  arrival \rl{or} departure route is assumed to be constant, \rl{but} the running time on different routes \rl{may differ}.}
\end{itemize}

\subsection{Model}
\label{subsec:mmodel}
The notation used to formulate the mathematical formulation is given in Tables~
\ref{notation_param}, and~\ref{notation_variable}.

\begin{table}[ht!]
\centering
\caption{Definition of sets and parameters.}
\label{notation_param}
\small{\small{
\begin{tabular}{p{69pt}p{366pt}}
\toprule
Notation & Definition\\
\midrule
$N$&Set of physical nodes, indexed by $i,j$\\
$N^a,N^d,N^{st}_s,N^{st}_m$ &Set of entering nodes, leaving nodes, siding nodes and mainline nodes\\
$F$& Set of trains, indexed by $f,f'$\\
$T$&Set of time periods, indexed by $t$\\
$A_f$& Set of macroscopic arcs that can be used by train $f$\\
$A_f^a,A_f^d,A_f^{st}$&Set of arrival, departure and siding waiting  arcs that can be used by train $f$\\
$\delta_f^-(v),\delta_f^+(v)$&Set of  arcs that originate from or end at vertex $v$ for train $f$\\
$w_1,w_2$& Weights of the two objectives\\
$T_g$& Running time of arc $g$\\
$T^a_f,T^d_f$&Desired arrival and departure times of train $f$\\
$\underline{D}^f,\overline{D}^f$&Minimum and maximum dwell time for train $f$ \\

$\Delta s$&Minimum headway required for the usage of the same 
siding track\\
\bottomrule
\end{tabular}}}
\end{table}

\begin{table}[ht!]
\centering
\caption{Definition of variables.}
\label{notation_variable}
\small{\small{
\begin{tabular}{p{35pt}p{400pt}}
\toprule
Notation & Definition\\
\midrule
$x_g^f$& Whether train $f$ uses macroscopic space-time arc $g$ (=1 yes; =0 no)\\
$\theta_l^f$& Whether train $f$ uses space-time resource $l$ at the microscopic level (=1 yes; =0 no)\\
$y_l^f$& Whether train $f$ uses siding resource $l$ implicitly at the microscopic level (=1 yes; =0 no)\\
$\mu_l^f$& Whether train $f$ uses siding resource $l$ at the microscopic level due to the lock of the arrival route (=1 yes; =0 no)\\
\bottomrule
\end{tabular}}}
\end{table}

\subsubsection{Objective}
The objective function minimizes the weighted sum of the \qin{travel} times and \qin{shift} times of all involved trains. The 
\qin{shift} time includes the arrival and departure \qin{shift}s from the desired times in the ideal timetable. The objective function can therefore be stated as follows:

\begin{equation}
    \min Z_1=w_1\sum_{f\in F}{\sum_{g\in A_f} {T_g \cdot x_g^f}}+w_2\left(\sum_{f\in F}|{\sum_{g\in A_f^a}{\hat{d_g}\cdot x_g^f}-T_f^a}|+\sum_{f\in F}|{\sum_{g\in A_f^d}{\hat{o_g}\cdot x_g^f}-T_f^d}|\right)\label{eq:Z1_obj}
\end{equation}

\subsubsection{Constraints}

Any feasible solution to the model must satisfy the following constraints:

\vspace{6pt}
\noindent$\bullet$ {\bf Train flow balance constraints at the macroscopic level}
\vspace{2pt}

\begin{equation}
    \sum_{g\in \delta_f^-(v)}{x_g^f} -\sum\limits_{g\in \delta_f^+(v)}{x_g^f}=
    \begin{cases}
    1&v=\overline{o} \\
    -1&v=\overline{k} \\
    0&\qin{v\in V_f\backslash(\overline{o}\cup \overline{k})}
    \label{train_flow}
    \end{cases}
    , \forall f \in F
\end{equation}

Constraints~\eqref{train_flow} ensure flow balance at the source, intermediate, and sink nodes. Only one available path is assigned to train $f\in F$ in the macroscopic space-time level. Unlike \citet{wang2023joint}, we do not impose such constraints on the utilization of microscopic space-time resources.

\vspace{6pt}
\noindent$\bullet$ {\bf Train dwell time constraints at the macroscopic level}

\begin{align}
    \sum_{g\in  A_f^{st}}{x_g^f} &\ge (1-\qin{x_{g_p}^f})\cdot \underline{D}^f \quad \quad  &\forall f\in F\label{min_dwell_constraint}\\
\sum_{g\in  A_f^{st}}{x_g^f} &\le (1-\qin{x_{g_p}^f})\cdot \overline{D}^f & \forall f\in F\label{max_dwell_constraint}
\end{align}

Constraints~\eqref{min_dwell_constraint} ensure the minimum dwell time of train $f\in F$ at the station, allowing sufficient time for passengers to board or alight the train. Constraints~\eqref{max_dwell_constraint} define the maximum allowable dwell time for train $f\in F$ at the station. If none of the siding waiting arcs can be occupied by $f\in F$, a virtual path is assigned to this train.

\vspace{6pt}
\noindent$\bullet$ {\bf Linking constraints between the two levels}

\begin{align}
    \qin{|\phi^{SG}_g|}\cdot x_g^f&\le \sum_{\qin{l\in \phi^{SG}_g}}{{\theta_l^f}}& \forall g\in A^a\cup A^d, \forall f\in F\label{map_uproute_lowsg}\\
    \qin{|\phi^{st}_g|}\cdot{x_g^f}&\le \sum_{\qin{l\in \phi^{st}_g}}{\mu_{l}^{f}}& \forall g\in A_f^a|d_g\in N^{st}_{s},\forall f\in F\label{map_between_arrival_siding}\\
    y_l^f&=\sum_{\qin{g\in A_l^d}}{x_{g'}^f} &\forall l\in B_f^{st}, \forall f\in F \label{map_between_o_i}\\
    \qin{|\phi_g^{ss}|\cdot x_g^f}&\qin{= \sum_{l\in \phi_g^{ss}}{\theta_l^f}}& \qin{\forall g\in A_f^{st}, \forall f\in F}\label{map_between_macro_micro_siding}
\end{align}

When a train occupies an inbound route or outbound route, the minimum headway $\Delta r$ required for the following train to use the same \gls{SG} resource in the bottleneck area is defined as the sum of the clearing time of the resource used by current train and the formulation time of the following train's route. As shown in Fig.~\ref{fig:map_two_ST_net}, the minimum headway between two consecutive trains can be considered as the implicit occupation of the \gls{SG} resources by the first train. Constraints~\eqref{map_uproute_lowsg} define the connection between the occupation of arrival (departure) arcs at the macroscopic level and the duration of occupying corresponding microscopic \gls{SG} resources at the microscopic level, which effectively performs the interlocking configuration. 

Constraints~\eqref{map_between_arrival_siding}- \eqref{map_between_macro_micro_siding} are the linking constraints used to connect the decision variables associated with the occupation of macroscopic arcs and microscopic siding resources. As shown in Fig.~\ref{fig:occupation_siding}, the siding arc is locked as soon as the train starts to occupy the corresponding inbound route. Constraints~\eqref{map_between_arrival_siding} define the connection between the occupation of a macroscopic arrival arc and the corresponding microscopic siding resources. A safety headway is required for two consecutive occupations of the siding track. This headway can be modelled as the implicit occupation of the siding track by the former train after its departure, as described in~\citet{zhang2020simultaneously}. \qin{We define the set of departure arcs that implicitly occupy the microscopic siding resource $l\in B^{st}$ as $A^d_l=\left\{g|g\in A_f^d,o_{g}=o_l,\hat{o}_{g}\in (\hat{o}_l- \Delta s,\hat{o}_l]\right\}$.} Constraints~\eqref{map_between_o_i} represent the connection between the occupation of a departure arc at the macroscopic level and the implicit occupation of the siding resources at the microscopic level. Constraints~\eqref{map_between_macro_micro_siding} state the connection between the actual occupations of siding tracks at the two levels \qin{where $\phi_g^{ss}=\left\{l|l\in B^{st},n_l=o_g,\hat{o}_l\in [\hat{o}_g,\hat{d}_g)\right\}$ represents the microscopic siding resources that correspond to the macroscopic siding arc $g\in A^{st}$}. \qin{The difference in time granularity between the macroscopic and microscopic levels is reflected in the construction of $\phi^{SG}_g$, $\phi^{st}_g$, $A_l^d$ and $\phi^{ss}_g$ and these linking constraints work under any granularity condition.} 
Clearly, the proposed two-level space-time network facilitates a straightforward and accurate computation of the occupation of track resources by a train given the interlocking system used at a complex station. 

\vspace{6pt}
\noindent$\bullet$ {\bf Capacity constraints for the track resources at the microscopic level}

\begin{align}
    \sum_{f\in F|l\in B_f^s}{\theta_l^f}\le 1\quad \quad \forall l\in B^s\label{headway_route}\\
    \sum_{f\qin{\in F|}l \in B^{st}_f} {(\theta^f_l  + y^f_l+\mu_{l}^{f})}  \le 1 \quad \quad \forall l \in B^{st}\label{siding_capacity}
\end{align}

To ensure the safety of train movements within the station, it is imperative that all microscopic station resources can only be occupied by a single train at any time instant.  Constraints~\eqref{headway_route} guarantee that only one train can occupy a particular \gls{SG} resource at any time instant at the microscopic space-time level. Constraints~\eqref{siding_capacity} ensure that the total number of all types of occupations of one microscopic siding resource at any time instant does not exceed one.





\vspace{6pt}
\noindent$\bullet$ {\bf Variable domains constraints}

Finally, the domains of all decision variables are given below.
\begin{align}
x^f_g & \in \left\{ {0,1} \right\} & \forall f\in F, \forall g \in A_f \label{arcdomain} \\
\theta^f_l &\in \left\{ {0,1} \right\} &\forall f \in F, \forall l \in B^{s}_f \cup B^{st}_f\label{resource_micro}\\
y^f_l  &\in \left\{ {0,1} \right\} & \forall f \in F, \forall l \in B^{st}_f \label{implicitdomain} \\
\mu^f_l &\in \left\{ {0,1} \right\} &\forall f \in F, \forall l \in B^{st}_f \label{a_siding_map}
\end{align}

\subsubsection{Model Discussion}
\label{balance_discussion}
A balanced resource utilization 
plays a significant role in traffic management because the use of platform tracks determines how frequently connected switches will be used. Specifically, more maintenance and replacement operations are required for switches that are overused. In contrast, a failure indicator for utilization will be triggered if a switch is rarely used. Our proposed model can be extended to consider the balanced use of platform tracks. There are two ways to realize this extension of even track use:

(1) Balanced track use can be encouraged by the objective function. For example, we can minimize the standard deviation or absolute value of the differences in occupation times between the platform tracks. Unfortunately, this approach can make the model difficult to solve. 

(2) Additional constraints can be enforced to incorporate the consideration of balanced track occupation. For example, we can translate balanced track use into the maximum number of trains allocated to a certain siding track. To achieve this, the maximum number of train occupations on one siding track should be not more than than the sum of the average usage $\bar{n}_{track}$ and the tolerance of imbalance $\bar{\theta}_{track}$, as shown in constraints~\eqref{balance_constraint}. This set of constraints can be adapted according to the preferences of railway planners.

\begin{equation}
    \sum_{f\qin{\in F|}g\in A_f}{\sum_{g\in A_f^{a}|d_g=i}{x_g^f}}\le \bar{n}_{track}+\bar{\theta}_{track}\quad \forall i\in N^{st}_s \label{balance_constraint}
\end{equation}

\section{Algorithm}
\label{sec:alg}
\subsection{Linearization}
To linearize the objective function in our proposed model, we introduce two positive integer variables, $y_f^{\text{shift,a}}$ and $y_f^{\text{shift,d}}$, for train $f\in F$. The objective function in the \gls{TLSTN} model can be rewritten as follows. The proof is not provided in this paper as it can be easily verified.

\begin{equation}
    \min Z_2=w_1\sum_{f\in F}{\sum_{g\in A_f} {T_g \cdot x_g^f}}+w_2\left(\sum_{f\in F}{y_{f}^{\text{shift},a}}+\sum_{f\in F}{y_{f}^{\text{shift},d}}\right)\label{ub_objective}
\end{equation}
where
\begin{align}
     y_{f}^{\text{shift},a}&\ge \sum_{g\in A_f^a}{\hat{d_g}\cdot x_g^f}-T_f^a&\quad \quad\forall f \in F\label{shift_a_1}\\
     y_{f}^{\text{shift},a}&\ge T_f^a-\sum_{g\in A_f^a}{\hat{d_g}\cdot x_g^f}&\quad \quad\forall f \in F\label{shift_a_2}\\
     y_{f}^{\text{shift},d}&\ge \sum_{g\in A_f^d}{\hat{o_g}\cdot x_g^f}-T_f^d&\quad \quad\forall f \in F\label{shift_d_1}\\
     y_{f}^{\text{shift},d}&\ge T_f^d-\sum_{g\in A_f^d}{\hat{o_g}\cdot x_g^f}&\quad \quad\forall f \in F\label{shift_d_2}\\
    y_{f}^{\text{shift},d}&\in \qin{Z} &\forall f \in F\label{shift_d_domain}\\
     y_{f}^{\text{shift},a}&\in \qin{Z} &\forall f \in F\label{shift_a_domain}
\end{align}

\subsection{Decomposition method}
\subsubsection{\qin{Solution framework}} 
\label{subsubsec:solution_framework}
\gls{BD} and \gls{LR} are two classical decomposition algorithms for solving large-scale optimization problems. \gls{BD}, which is also called variable partitioning, was initially proposed to tackle 
linear programming problems, but is also effective at solving mixed-integer linear programmings with a so-called {\it dual block angular structure.} 
Variables in problems with such a structure can typically be split into two sets: variables $X$ and variables $Y$. This enables the original problem to be decomposed into a master problem over $X$ and a subproblem over $Y$ given the fixed $X$ values. The subproblem returns information to the master problem by way of optimality and feasibility cuts, and the two problems are iteratively solved until the suproblem can find no violated inequalities. \gls{LBBD} is an extension of \gls{BD} and allows for more general subproblems using a logic formalism (e.g., with discrete decision variables)~\citep{hooker1995testing,hooker2003logic}. 
The 
\gls{LBBD} method 
has been applied to solve railway scheduling and rescheduling problems (e.g.,~\citet{keita2020three,leutwiler2022logic,lamorgese2015exact,lamorgese2016optimal}). However, \gls{LBBD} provides no standard scheme for generating Benders cuts, meaning logic-based cuts are generally problem specific~\citep{hooker2007planning}. Furthermore, the dynamically added cuts typically destroy the decomposability of the master problem. 

In contrast,  
\gls{LR} algorithm is a common decomposition algorithm for solving large-scale integer programming problems. This algorithm dualizes the so-called hard constraints of a model to the objective function through Lagrange multipliers and further decomposes the problem into a series of independent subproblems. Any violations of hard constraints are reflected by Lagrange multipliers. The greater the magnitude of violation for a specific constraint, the larger the penalty in the objective function. This approach has been widely applied to optimize railway operations (e.g.,~\citet{zhang2020joint},\citet{luan2017integrated}, \citet{xu2021train}, \citet{zhang2021collaborative},~\citet{meng2014simultaneous}). In these studies, the capacity or headway constraints between trains are moved into the objective function, and then the relaxed model can be split into a set of train-specific resource-constrained shortest-path subproblems. This approach maintains the decomposability of problems and facilitates the application of many well-proven algorithms. 


Formally, we introduce two sets of Lagrange multipliers $\lambda_l, \forall l\in B^s$ and $\lambda_l, \forall l\in B^{st}$ to dualize \qin{microscopic-exclusive capacity} constraints~\eqref{headway_route} and~\eqref{siding_capacity}, respectively, into the objective function and penalize their violation. The following model is referred to as the \gls{TLSTN}-LR model.

\begin{equation}
    \min \qin{Z_{LR}}=Z_2+\sum_{l\in B^s}{\lambda_l\left(\sum_{f\in F|l\in B_f}{\theta_l^f}-1\right)}+\sum_{l \in B^{st}}{\lambda_l\left(\sum_{f\qin{\in F|}l \in {B_f^{st}}} {(\theta^f_l  + y^f_l+\mu_{l}^{f})}  - 1\right)} \label{z_lr}
\end{equation}

\noindent subject to constraints~\eqref{train_flow}-\eqref{map_between_macro_micro_siding},~\eqref{arcdomain}-\eqref{a_siding_map},~\eqref{shift_a_1}-\eqref{shift_a_domain} and~\eqref{lambda_l}-\eqref{lambda_g}. The domains of the Lagrange multipliers are defined as:

\begin{align}
    \lambda_l&\ge 0, \forall l\in B^{s}\label{lambda_l}\\
    \lambda_l&\ge 0, \forall l\in B^{st}\label{lambda_g}
\end{align}

According to constraints~\eqref{map_uproute_lowsg}, for train $f\in F$, if $\theta^f_l=1$, then only one macroscopic arc $g\in A^{a}_f\cup A^{d}_f$ that occupies microscopic resource $l$ is selected by train $f$ (i.e., $x_g^f=1$). Therefore, the second part of the objective function in~\eqref{z_lr} can be rewritten as follows:

\begin{align}
    \sum_{l\in B^s}{\lambda_l\left(\sum_{f\in F|l\in B_f^s}{\theta_l^f}-1\right)}&=\sum_{f\in F}{\sum_{l\in B_f^s}{\theta_l^f{\lambda_l}}}-\sum_{l\in B^s}{\lambda_l}\nonumber\\
    &= \sum_{f\in F}{\sum_{g\in A_f^{a}\cup A_f^{d}}{x_g^f\left(\sum_{l\in \phi_g^{SG}}{\lambda_l}\right)}}-\sum_{l\in B^s}{\lambda_l}\nonumber\\
    &= \sum_{f\in F}{\sum_{g\in A_f^{a}\cup A_f^{d}}{\widetilde{c}_g{x_g^f}}}-\sum_{l\in B^s}{\lambda_l}
\end{align}

\qin{It should be noted that in constraints~\eqref{map_between_macro_micro_siding}, $\sum_{l\in \phi_g^{ss}}{\theta _l^f}$ is \rl{essentially} a copy of the $x_g^f$ variable for train $f$ and arc $g\in A^{st}_f$.}
Finally, 
the objective function~\eqref{z_lr} can be written in an aggregated form as~\eqref{Z_LR_f}. 

\begin{align}
    \min \qin{Z_{LR}}&=Z_2+ \sum_{f\in F}{\sum_{g\in A_f^{a}\cup A_f^{d}}{\widetilde{c}_g{x_g^f}}}-\sum_{l\in B^s}{\lambda_l}+\sum_{l \in B^{st}}{\lambda_l\left(\sum_{f\qin{\in F|}l \in {B_f^{st}}} {(\theta^f_l  + y^f_l+\mu_{l}^{f})}  - 1\right)} \nonumber\\
    &=\sum_{f\in F}\left\{{\sum_{g\in A_f} {\hat{c}_g^f \cdot x_g^f}}\right\}-\sum_{l\in B^s\cup B^{st}}{\lambda_l}\nonumber\\
    &=\sum_{f\in F}{\left\{\qin{Z_{LR}^f}\right\}}-\sum_{l\in B^s\cup B^{st}}{\lambda_l}\label{Z_LR_f}
\end{align}

\noindent where 
\begin{equation}
    \hat{c}_g^f=
    \begin{cases}
    w_1 T_g+\widetilde{c}_g+w_2|\hat{d}_g-T^a_f|+\sum_{l\in B^{st}_f|o_l=d_g,\hat{o}_l\in [\hat{o}_g,\hat{d}_g)}{\lambda_l}& \forall g \in A_f^{a},\\
    w_1 T_g+\widetilde{c}_g+w_2|\hat{o}_g-T^d_f|+\sum_{l\in B^{st}_f|o_l=o_g,\hat{o}_l\in [\hat{o}_g,\hat{o}_g+\Delta_s)}{\lambda_l}& \forall g \in A_f^{d},\\
    w_1 T_g+\qin{\sum_{l\in \phi_g^{ss}}{\lambda_l}}& \forall g \in A_f^{st},
    \end{cases}
    \qquad \forall f \in F\label{cost_arcs}
\end{equation}

The \qin{macroscopic-exclusive} problem
\qin{of the \gls{TLSTN}-LR model} is composed of a set of train blocks, each of which can  be solved independently. Each train-block entails solving a resource-constrained shortest-path problem, which can be effectively handled by dynamic programming given the values of the Lagrange multipliers. For the sake of completeness, we present the mathematical model for block $f$ \qin{below}. 

 \begin{align}
 \min & \quad \qin{Z_{LR}^f} \nonumber\\
 \text{s.t.}&
\begin{cases}
\text{Constraints}~\eqref{train_flow} -~\eqref{max_dwell_constraint}\\
\text{Constraints} ~\eqref{arcdomain}\\
\end{cases}
\label{ITSMSODM}
\end{align}

\begin{figure}[h!]
    \centering
    \includegraphics[scale=0.6]{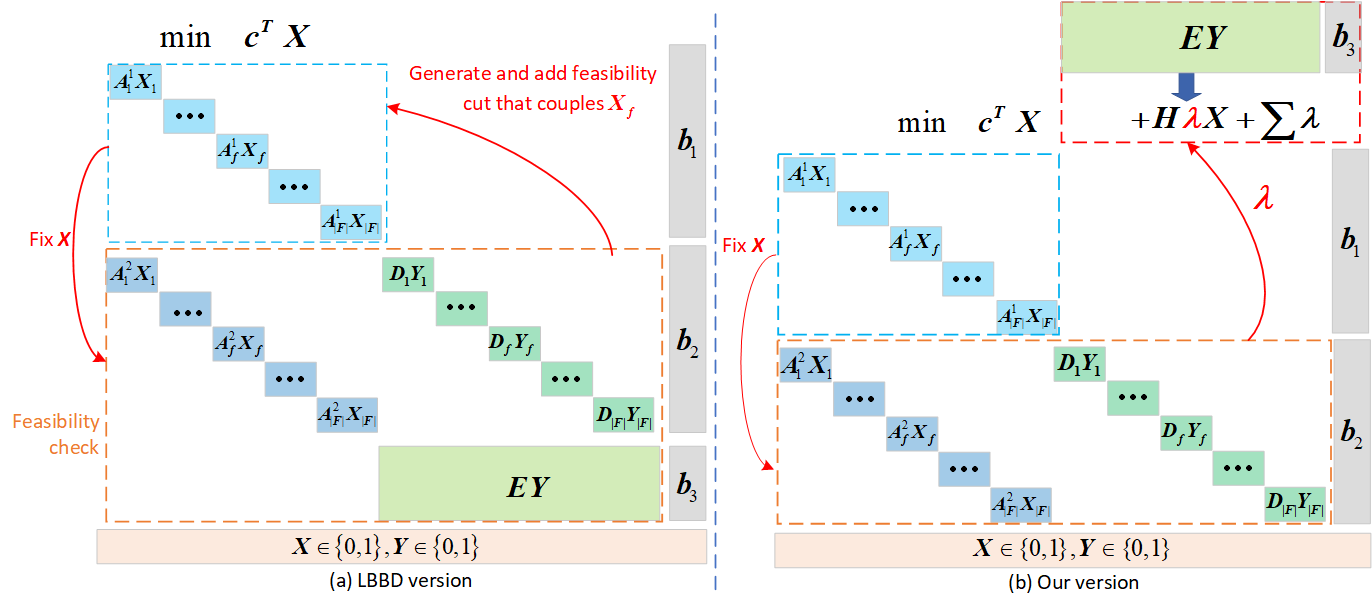}
    \caption{An overview of the \gls{LBBD} and \gls{2-L LR} decomposition approaches.}
    \label{fig:model_structure}
\end{figure}

In Fig.~\ref{fig:model_structure}, we categorize the constraints of the \gls{TLSTN} model into three groups. The first group are macroscopic-exclusive constraints that enforce path selection for trains without considering any microscopic capacity headway. The second group comprises linking constraints, which transform the macroscopic arc choice into the occupation of microscopic resources. The constraints in these two groups are train independent. The final group of constraints is microscopic-exclusive and ensures capacity requirements at the microscopic level are respected. These constraints couple different trains and make the model difficult to solve.

We illustrate a \gls{LBBD} approach for the proposed \qin{\gls{TLSTN}} model in Fig.~\ref{fig:model_structure}-a). The master problem includes the objective function and the constraints in the first group. The subproblem, however, is a feasibility problem, which includes the fixed solution of the master problem and the remaining constraints for ensuring the feasibility of the master solution at the microscopic level. If the path selection is infeasible, then a feasibility cut in the form of a no-good cut \rl{is} generated and added to the master problem. The master problem is solved iteratively until the result is feasible. \rl{We use no-good cuts of the form given in~\eqref{eq:ngcut}. These are typically used in applications of \gls{LBBD}, see e.g., \cite{Hooker2019, forbes2024840}, and can be used to eliminate feasible solutions to the master problem that lead to subproblem infeasibility.}

\begin{equation}
\label{eq:ngcut}
    \sum_{ \qin{x\in \hat{X}}} \qin{x}\leq |\qin{\hat{X}}|-1
\end{equation}
where $\qin{\hat{X}=\left\{x_g^f|x_g^f=1,\forall g\in A_f, \forall f\in F\right\}}$. \rl{This set contains the decision variables that are associated with the arcs that are used in the macroscopic-level-feasible solution. Since this solution induces infeasibility for the microscopic level, \eqref{eq:ngcut} is used to eliminate precisely this solution from the master problem.} 
While \rl{such} cuts will ultimately ensure an optimal solution to the non-decomposed problem is obtained (by iteratively cutting off the \rl{solutions that induce microscopic-level infeasibility)}, no information on where the infeasibility is provided by the cut, and the decomposability of the master problem further into train independent problems is \rl{lost}.

Fig.~\ref{fig:model_structure}-b), on the other hand, illustrates \qin{our proposed} \gls{2-L LR} approach. We propose to dualize the capacity constraints at the microscopic level and penalize their violation in the objective function using Lagrange multipliers. The 
\qin{\gls{TLSTN}-LR} model therefore contains the macroscopic level problem and the linking constraints only. We are left with a decomposition that resembles the \gls{LBBD} decomposition in that the relaxed problem can be decomposed into a macroscopic level train pathing master problem, and a subproblem that links this solution to the microscopic level. We note here that for a fixed master solution the subproblem admits a unique solution and feasibility is not an issue. Instead of returning a feasibility cut to the master problem, we can penalize the use of overutilized microscopic level track resources by adjusting the corresponding Lagrange multipliers. Such an approach has two main advantages. Firstly, with this approach we can guide the train pathing master problem by indicating where the microscopic level infeasibility arises from. Secondly, by not adding cuts, we retain the decomposability of the master problem into train independent problems.
Finally, we note that the Lagrange multipliers are each associated with an individual microscopic level track resource. Since a macroscopic level element typically contains multiple microscopic level resources, the penalty associated with the use of a macroscopic level track resource is the sum of the Lagrange multipliers associated with corresponding microscopic level resources.

\begin{figure}[h!]
    \centering
    \includegraphics[scale=0.91]{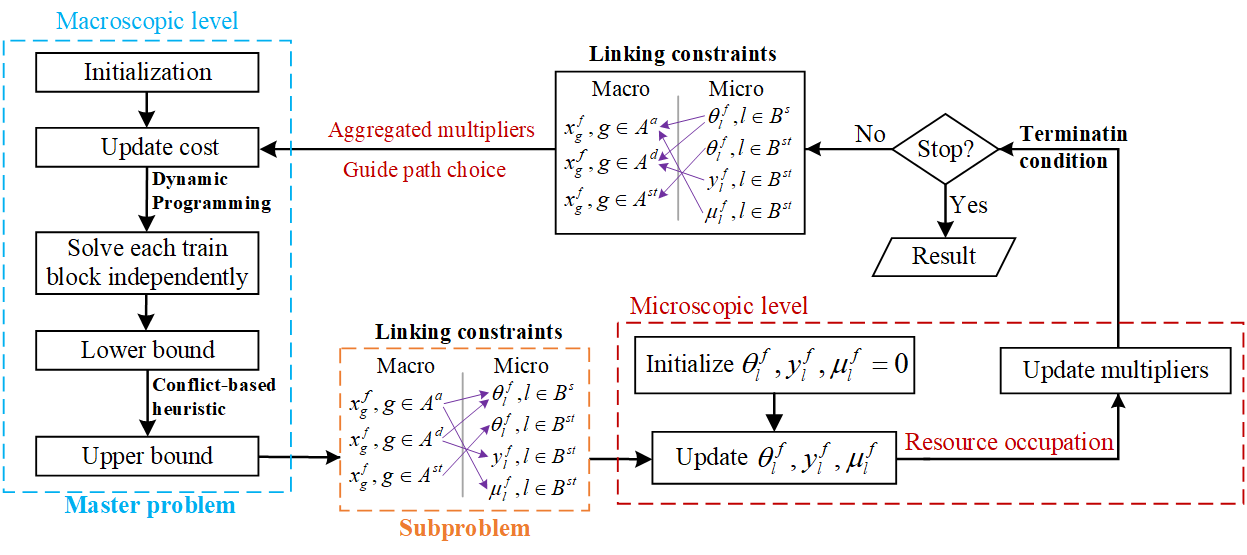}
    \caption{Solution framework of our proposed \gls{2-L LR}.}
\label{fig:two_layered_LR}
\end{figure}

An overview of the proposed \gls{2-L LR} algorithm is shown in Fig.~\ref{fig:two_layered_LR}. The subproblem transforms the occupations of macroscopic arcs in the master solution into occupations of microscopic resources. The multipliers are then updated accordingly based on the violation of capacity constraints at the microscopic level.
Next, the multipliers are aggregated into the costs of the corresponding macroscopic arcs to guide the path choice of trains in the master problem. A more detailed description of this approach is provided in~Algorithm~\ref{alg:framework}. A sub-gradient method is used to update the Lagrange multipliers, and the step size retains the same value after iteration $m_\alpha$. Readers can 
\qin{refer to \citet{held1974validation}, \citet{boyd2003subgradient} and \citet{fisher2004lagrangian}} for additional details on this method. In addition, we apply a  \qin{priority-rule based }heuristic method to generate a feasible upper bound. 
\qin{We assign \rl{a} higher priority to trains \rl{that have fewer} conflicts with other trains in the lower bound solution and then sequentially schedule the trains based on \rl{their} priority.} When trains with higher priority are scheduled, \rl{macroscopic arcs that conflict with the resources used by the scheduled trains, as determined by constraints~\eqref{headway_route} and~\eqref{siding_capacity}, are removed} from the space-time network, and subsequent trains can only occupy the remaining available resources \qin{that \rl{do} not conflict with the previously scheduled trains. \rl{This ensures a feasible solution is obtained}~\citep{zhang2020simultaneously}}. 

\begin{small}
    \begin{algorithm}[ht!]
\caption{Two-level Lagrangian Relaxation} \label{alg:framework}
\begin{algorithmic}[1]
\State \textbf{//Step 1: Initialization}
  \State  \quad Initialize iteration number $m=0$,  multipliers $\lambda_l^m=0, \forall l\in B^{s}\cup B^{st}$, step size $\alpha(m)=1$;
  \State \quad Set best upper bound ${\text{UB}}^*= +\infty$ and corresponding solution $\left\{X_{\text{UB}}^*\gets \O\right\}$,  best lower bound ${\text{LB}}^*= -\infty$ and corresponding solution $\left\{X_{\text{LB}}^*\gets \O\right\}$.
\State \textbf{//Step 2: Generate lower bound solution}
  \State \textbf{\quad //Step 2.1: Update cost of arcs at the macroscopic level}
    \State{\quad\quad \textbf{for} $g \in A$ \textbf{do}}
        \State{\quad\quad\quad Update aggregated cost $\hat{c}_g$ by Eq.~\eqref{cost_arcs}}.
  \State \textbf{\quad //Step 2.2: Solve each train block}
    \State \quad {Initialize current lower bound solution $\left\{X_{\text{LB}}\gets \O\right\}$}
    \State{\quad\quad \textbf{for} $f \in F$ \textbf{do}}
       \State{\quad\quad\quad $\left\{X_{\text{LB}}^{f}\gets \O\right\} $};
      \State{\quad\quad\quad Generate $\left\{X_{\text{LB}}^{f}\right\}$ based on dynamic programming method as described in~\citet{zhang2020simultaneously}};
      \State{\quad\quad\quad $\left\{X_{\text{LB}}\right\}=\left\{X_{\text{LB}}\right\}\cup \left\{X_{\text{LB}}^{f}\right\}$ }.
  \State \textbf{\quad //Step 2.3: Generate lower bound and corresponding solution}
    \State \quad\quad Calculate current lower bound $\text{LB}^m$ by Eq.~\eqref{Z_LR_f};
    \State \quad\quad \textbf{if} $\text{LB}^m>\text{LB}^*$ \textbf{then}
      \State \quad\quad\quad $\text{LB}^*= \text{LB}^m$, $\left\{X_{\text{LB}}^*\right\}= \left\{X_{\text{LB}}\right\}$.
\State \textbf{//Step 3: Generate upper bound solution}
\State \quad {Initialize current upper bound solution $\left\{X_{\text{UB}}\gets \O\right\}$}
  \State \quad \textbf{if} there is conflicting occupation in $\left\{X_{\text{LB}}\right\}$ \textbf{then}
    \State \quad\quad Generate $\left\{X_{\text{UB}}\right\}$ based on heuristic method;
  \State \quad \textbf{Else} 
  \State \quad\quad $\left\{X_{\text{UB}}\right\}=\left\{X_{\text{LB}}\right\}$;
    \State \quad Calculate current upper bound $\text{UB}^m$ by Eq.~\eqref{ub_objective};
    \State \quad \textbf{if} $\text{UB}^m<\text{UB}^*$ \textbf{then}
       \State \quad\quad $\text{UB}^*= \text{UB}^m$, $\left\{X_{\text{UB}}^*\right\}= \left\{X_{\text{UB}}\right\}$.
\State \textbf{//Step 4: Update Lagrange multipliers at the microscopic level}
\State \quad Calculate $\theta_l^f$ by Eq.~\eqref{map_uproute_lowsg} and~\eqref{map_between_macro_micro_siding};
\State \quad Calculate $\mu_l^f$ and $y_l^f$ by Eq.~\eqref{map_between_arrival_siding} and~\eqref{map_between_o_i}, respectively;
  \State \quad \textbf{for} $l \in B^{s}$ \textbf{do}
      \State \quad\quad $\lambda_l^{m+1}=\max\left\{0,\lambda_l^m+\alpha(m)\cdot \left(\sum_{f\in F|l\in B_f^s}{\theta_l^f-1}\right)\right\}$;
  \State \quad \textbf{for} $l \in B^{st}$ \textbf{do}
      \State \quad\quad $\lambda_l^{m+1}=\max\left\{0,\lambda_l^m+\alpha(m)\cdot \left(\sum_{f\in F|l \in {B_f^{st}}} {(\theta^f_l  + y^f_l+\mu_{l}^{f})}  - 1\right)\right\}$;
  \State \quad \textbf{if} $m$ is less than a certain iteration number $m_\alpha$ \textbf{then}
      \State \quad\quad $\alpha(m)=1/m+1$;
\State \textbf{//Step 5: Evaluate solution quality and termination condition test}
    \State \quad \textbf{if} $\left\{X_{LB}\right\}$ is feasible and \textbf{$\lambda$}=\textbf{0} \textbf{then}
    \State \quad \quad ${\text{UB}}^*={\text{LB}}^*$, $\left\{X_{\text{UB}}^*\right\}= \left\{X_{\text{LB}}^*\right\}$;

    \State \quad $\text{GAP}^{m}=\frac{\text{UB}^*-\text{LB}^*}{\text{UB}^*}$
    \State \quad If $\text{GAP}^{m}$ is sufficiently small or $\text{UB}^*-\text{LB}^*<1$ or 
    $m$ reaches 
    max iterations, or the computing time exceeds the 
    limit, stop; Otherwise, $m=m+1$, go back to Step 2.
\end{algorithmic}
\end{algorithm}
\end{small}

\subsubsection{Dynamic multiplier generation}
The number of Lagrange multipliers is large. The values of all multipliers are initially set to zero and then adapted according to the train occupation in each iteration. If any occupation of a specific resource violates the capacity constraint in the current iteration, then its corresponding Lagrange multipliers will be increased and the resource will be used at a higher cost in the next iteration. However, if the utilization of a certain resource satisfies the microscopic level capacity requirements, then its corresponding Lagrange multipliers will be kept at the current values or reduced to a non-negative value~\citep{zhang2020simultaneously}. Explicitly handling all multiplier updates at  every iteration is tedious and will result in poor algorithm performance.

At every iteration, most multipliers have  the value zero. Therefore, we adopt a dynamic multiplier approach and only generate multipliers as needed to reduce the total number of multipliers~\citep{zhang2020joint,caprara2002modeling} being explicitly considered. We introduce a resource pool $M^m$ to define the space-time resources with Lagrange multipliers greater than zero at the microscopic level in iteration $m$. Resource pools $M^0$ are initially set to empty; i.e., we implicitly set all multipliers to zero. Then, the resource pools $M^m$ are updated based on the occupations of resources by trains after solving all train blocks in each iteration. The method for dynamically updating space-time resources with multipliers greater than zero is presented in~Algorithm~\ref{alg:dynamic_multiplier} in~\ref{sec:appendix_dp}.

\subsubsection{Consideration of balanced occupation}
If we consider balanced track use by incorporating constraints~\eqref{balance_constraint}, then we should also dualize these constraints into the objective function as they impact/couple sets of train paths. We, therefore, introduce multipliers $\lambda_i\ge0,\forall i \in N_s^{st}$ to penalize the corresponding constraint violation in the objective function. Because the number of balance constraints is equal to the number of siding tracks, we explicitly consider all multipliers in this set in every iteration. The specific modifications that must be made to the proposed \gls{2-L LR} framework are provided in~\ref{sec:appendix_balance}.

\section{Computational experiments}
\label{sec:experiment}
 We evaluate 
the performance of our proposed two-level model and the \gls{2-L LR} algorithm for several cases based on two stations of different sizes. The first station is a small virtual station, while the second is an actual railway station located on one of the busiest Chinese high-speed railway lines. The algorithm was implemented in Python 3.7 using Spyder and executed on a personal computer with an Intel (R) Core (TM) i7-1165G7 @ 2.80 GHz processor and 16 GB of memory. Additionally, we compare the performance of our methodology with several other methods. The version of CPLEX used in our experiments is 12.9.0.

\subsection{Tests on a small virtual station}
\begin{figure}[h!]
    \centering
    \includegraphics[scale=0.97]{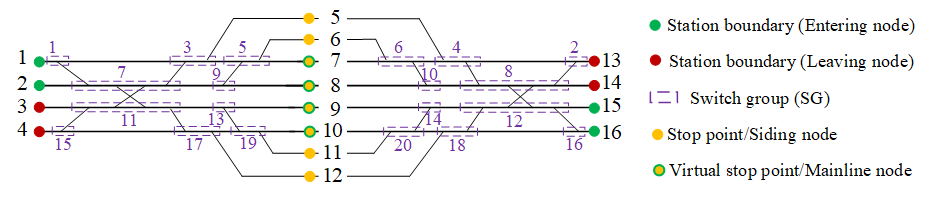}
    \caption{Station layout and macroscopic physical network in small case.}
    \label{fig:station_small_case}
\end{figure}

We first consider the small artificial station shown in Fig.~\ref{fig:station_small_case}. 
This station connects two railway lines and contains four siding tracks and four mainline tracks, respectively. The \gls{SG} resources are also given and are marked with odd numbers for the left side of the station and even numbers for the right side. All trains traverse the station from an entering node to a leaving node. Trains can only dwell on siding tracks. A flexible track utilization rule is applied in this study; stopping trains from one direction can dwell on any siding tracks, whereas the mainline tracks can only serve nonstop trains from one end of the mainline to the other end. 
The macroscopic physical network constructed for this small virtual station consists of 20 inbound arcs and 20 outbound arcs. The planning time horizon is 40 minutes and 
the time discretization is set to 15 seconds. The minimum headway required to use the same resource is assumed to be 30 seconds. The values of the weights in the objective function are set to one.

To analyze the performance of the proposed \gls{2-L LR} algorithm, we investigate 17 cases and compare the results to those of four other methods. First, we obtain the results by applying the commercial solver CPLEX directly to the full problem. The second method is the \gls{LBBD} method described in Section~\ref{subsubsec:solution_framework}. The \gls{LBBD} method with no-good cuts is implemented using cut callbacks at nodes that yield integer solutions during the branch-and-cut process in CPLEX. The third method is similar to the second but implements an iterative version of \gls{LBBD}. In other words, the master problem is fully solved to integer optimality before assessing subproblem feasibility.
We also test a fourth method, which we term the Lazy method. In this approach, we build the non-decomposed problem, as in method (1), however, we specify the linking constraints as lazy constraints in CPLEX. Performance comparisons are presented in Table~\ref{performance_comparison_small_case}. The first column shows the case ID. For each case, additional information such as the number of trains, space-time paths, 
conflicts between each pair of train paths, the number of variables and constraints, is presented in Table~\ref{small_case_details}. We evaluate the performance of the solutions on instances with the same number of trains but different characteristics in terms of origins, destinations, and arrival and departure distributions. For each method, we report the best upper bound (UB), the best lower bound (LB), the gap between UB and LB ($\Delta$), and the solution time ($t$) in seconds. For the 
\gls{2-L LR} method, we also present the total number of iterations (Iter) and the time at which UB is found ($t_u$). The termination condition we used in our proposed approach is a maximum of 1500 iterations whereas the condition for the other four methods is a maximum solving time of 1800 seconds. 

As can be observed, the performances of CPLEX and the Lazy method are similar; however, generally, the former outperforms the latter. When the number of trains is small, such as in Cases 1-11, both methods can obtain an optimal solution within the time limit. Solution times do vary for instances with the same number of trains, such as Cases 5 and 6. This is because the trains in these two cases have different characteristics, which can give rise to more or fewer paths and conflicting movements. The number of train paths and the number of conflicting movements in the latter case are both greater than those of the former. As the number of trains increases,  as in Cases 12-14, the optimality gap for CPLEX remains in the range of 35\% to 40\%, whereas that of the Lazy method increases from 57.6\% to 91\%. For the remaining cases, the gaps of CPLEX increase from 40.8\% to \qin{45.0\%}. In contrast, the Lazy method cannot obtain a feasible solution within the time limit. This is because an increase in problem size increases the number of linking constraints. The \gls{LBBD} and iterative approaches with no-good cuts perform much worse overall. These two methods can only provide upper bounds for 5 and 4 out of the 17 cases, respectively. Additionally, the lower bounds provided by these two methods are significantly weaker than those provided by the Lazy method.  This is because the number of linking constraints between levels is much less than the number of microscopic capacity constraints. In the Lazy approach, if we have to add every linking constraint, then we generate the full model. In contrast, we must solve the subproblem with capacity constraints in the other two 
methods, and feasibility cuts remove only one infeasible master solution at a time. Futhermore, no information on why the solution is infeasible is given to the master problem, meaning that the master problem may, on the next iteration, provide a structurally similar solution to the infeasible one. We also test the \gls{LBBD} with no-good cuts using callback at fractional nodes during the branch-and-cut process in CPLEX. The results are not much better than those obtained with callback at the integer nodes, so we do not discuss additional details here.

In Cases 1-11, where the numbers of trains are relatively low and CPLEX can obtain optimal solutions, our proposed algorithm is able to provide optimal solutions for nine cases. Although it is relatively weak in terms of proving optimality, the \gls{2-L LR} algorithm can provide optimal solutions or good upper bounds in a much shorter solving time. As the number of trains increases, the advantages of our designed algorithm gradually become more apparent. The upper bounds obtained by the proposed method are better than those obtained by CPLEX for large numbers of trains. Additionally, the latter requires significantly more solving time. As the number of trains increases, the gaps between the upper and lower bounds for CPLEX grow to 70.8\%. In contrast, the optimality gaps of the proposed approach remain between 11\% and 15\% as the size of the problem grows. Furthermore, the proposed method can provide a better upper and tighter lower bounds in a much shorter solving time than CPLEX in Cases 13-17. In particular, for Case 16, the \gls{2-L LR} algorithm provides a solution with an objective value of 6081 in 69 seconds, which is better than the value of 7318 obtained by CPLEX in 1800 seconds. Compared to the \gls{LBBD} approach with no-good cuts, our approach performs better in all small-sized cases.  

\renewcommand{\arraystretch}{1.5}
\newgeometry{top=1cm, bottom=1.8cm}
\begin{landscape}
    \centering
    \begin{table}
\centering
\caption{Performance comparison for the small-sized cases.}
\label{performance_comparison_small_case}
\small{
\small{
\small{
\small{
\begin{tabular}{r|rrrr|rrrr|rrrr|rrrr|rrrrrr}
\toprule
\multicolumn{1}{c}{Case}& \multicolumn{4}{c}{CPLEX}&\multicolumn{4}{c}{\gls{LBBD} \& no-good cuts}&\multicolumn{4}{c}{Iter \& no-good cuts}&\multicolumn{4}{c}{Lazy-linking}&\multicolumn{6}{c}{\gls{2-L LR}}\\
\cmidrule(lr){2-5} \cmidrule(lr){6-9} \cmidrule(lr){10-13}\cmidrule(lr){14-17}\cmidrule(lr){18-23} 
&UB&LB& $\Delta$ (\%) &$t$(s)&UB&LB& $\Delta$ (\%) & $t$ (s)&UB&LB&$\Delta$ (\%)& $t$ (s)&UB&LB& $\Delta$ (\%)& $t$(s)&UB&LB&$\Delta$  (\%)&$t$(s)&Iter&$t_u$(s)\\
\midrule
1&669&669&0.0&0.1&669&669&0.0&10&669&669&0.0&0.2&669&669&0.0&0.1&669&669&0.0&0.4&9&0.1\\
2&943&943&0.0&0.8&943&693&26.4&*&-&830&-&*&943&943&0.0&10&943&943&0.0&39&844&0.1\\
3&695&695&0.0&0.1&695&695&0.0&21&695&695&0.0&0.2&695&695&0.0&0.1&695&695&0.0&0.5&9&0.1\\
4&1460&1460&0.0&135&-&840&-&*&-&1033&-&*&1460&1460&0.0&457&1461&1457&0.3&75&$\dagger$&0.1\\
5&715&715&0.0&0.2&715&715&0.0&29&715&715&0.0&0.3&715&715&0.0&0.5&715&715&0.0&0.6&9&0.1\\
6&1626&1626&0.0&164&-&1004&-&*&-&1191&-&*&1626&1626&0.0&693&1626&1622&0.3&104&$\dagger$&2\\
7&741&741&0.0&0.2&741&741&0.0&41&741&741&0.0&0.3&741&741&0.0&0.1&741&741&0.0&0.6&9&0.1\\
8&1735&1735&0.0&183&-&1030&-&*&-&1217&-&*&1735&1735&0.0&1010&1735&1670&3.8&134&$\dagger$&10\\
9&1180&1180&0.0&8&-&911&-&*&-&1067&-&*&1180&1180&0.0&8&1180&1180&0.0&74&831&0.1\\
10&1755&1755&0.0&534&-&1050&-&*&-&1231&-&*&1755&1755&0.0&1465&1755&1706&2.8&144&$\dagger$&86\\
11&1796&1796&0.0&841&-&1076&-&*&-&1256&-&*&1796&1796&0.0&966&1803&1767&2.0&153&$\dagger$&34\\
12&3021&1952&35.4&*&-&1429&-&*&-&1611&-&*&3780&1603&57.6&*&3171&2764&12.9&244&$\dagger$&16\\
13&4067&2463&39.4&*&-&1903&-&*&-&2091&-&*&9446&2060&78.2&*&3832&3279&14.5&283&$\dagger$&122\\
14&4439&2884&35.0&*&-&2380&-&*&-&2577&-&*&28010&2533&91.0&*&4349&3782&13.1&350&$\dagger$&344\\
15&5809&3440&40.8&*&-&2872&-&*&-&3070&-&*&-&3048&-&*&5165&4534&12.2&438&$\dagger$&412\\
16&7318&4023&45.0&*&-&3360&-&*&-&3564&-&*&-&3523&-&*&6081&5353&12.0&477&$\dagger$&69\\
17&\qin{7710}&\qin{4362}&\qin{43.4}&*&-&3511&-&*&-&4010&-&*&-&3719&-&*&6367&5658&11.1&532&$\dagger$&83\\
\bottomrule
\multicolumn{8}{c}{Note: $*$ denotes 1800 seconds; $\dagger$ represents 1500 iterations.}
\end{tabular}}}}}
\end{table}
\end{landscape}
\renewcommand{\arraystretch}{1.5}
\restoregeometry

For the results of the \gls{2-L LR} approach in Table~\ref{performance_comparison_small_case}, an upper bound is generated at each iteration. We now compare this with a variant that only generates an upper bound in the final iteration. To increase the probability of obtaining a good upper bound, we randomly shuffled the order of trains with the same number of conflicting trains. Because the method used for finding the upper bound has no effect on the lower bound for a given iteration number, we only compare the upper bounds obtained by the two different upper bound generation strategies in Table~\ref{ub_comparison}. In the proposed method, for these small cases, generating an upper bound only in the final iteration can significantly reduce the total solving time required by the algorithm, by up to 47.1\%. However, the solution quality is reduced compared to the method of generating an upper bound in every iteration.

\renewcommand{\arraystretch}{1}
\begin{table}[ht!]
\centering
\caption{Performance comparison for the small-sized cases with different upper bound strategies.}
\label{ub_comparison}
\small{
\small{
\small{\small{
\begin{tabular}{rrrrrrr}
\toprule
\multicolumn{1}{c}{Case}& \multicolumn{2}{c}{UB-iterative}&\multicolumn{2}{c}{UB-final iteration}&\multicolumn{1}{c}{$\Delta$(\%)}&\multicolumn{1}{c}{$\Delta t$(\%)}\\
\cmidrule(lr){2-3} \cmidrule(lr){4-5} &UB1&$t1$(s)&UB2&$t2$(s)\\
\midrule
1&669&0.4&669&0.3&0&-25.0\\
2&943&39&943&42&0&7.7\\
3&695&0.5&695&0.4&0&-20.0\\
4&1461&75&1505&57&3.0&-24.0\\
5&715&0.6&715&0.5&0&-16.7\\
6&1626&104&1731&67&6.5&-35.6\\
7&741&0.6&741&0.5&0&-16.7\\
8&1735&134&1787&78&3&-41.8\\
9&1180&74&1180&61&0&-17.6\\
10&1755&144&1808&80&3.0&-44.4\\
11&1803&153&1900&106&5.4&-30.7\\
12&3171&244&3830&129&20.8&-47.1\\
13&3832&283&7342&159&91.6&-43.8\\
14&4349&350&7820&192&79.8&-45.1\\
15&5165&438&8771&238&69.8&-45.7\\
16&6081&477&15084&315&148.1&-34.0\\
17&6367&532&15422&346&142.2&-35.0\\
\bottomrule
\multicolumn{7}{l}{$\Delta=\frac{(\text{UB2}-\text{LB1)}}{\text{UB1}}*100$; $\Delta t=\frac{(t2-t1)}{t1}*100$}
\end{tabular}}}}}
\end{table}
\renewcommand{\arraystretch}{1.5}

\begin{figure}[h!]
    \centering
    \includegraphics[scale=0.49]{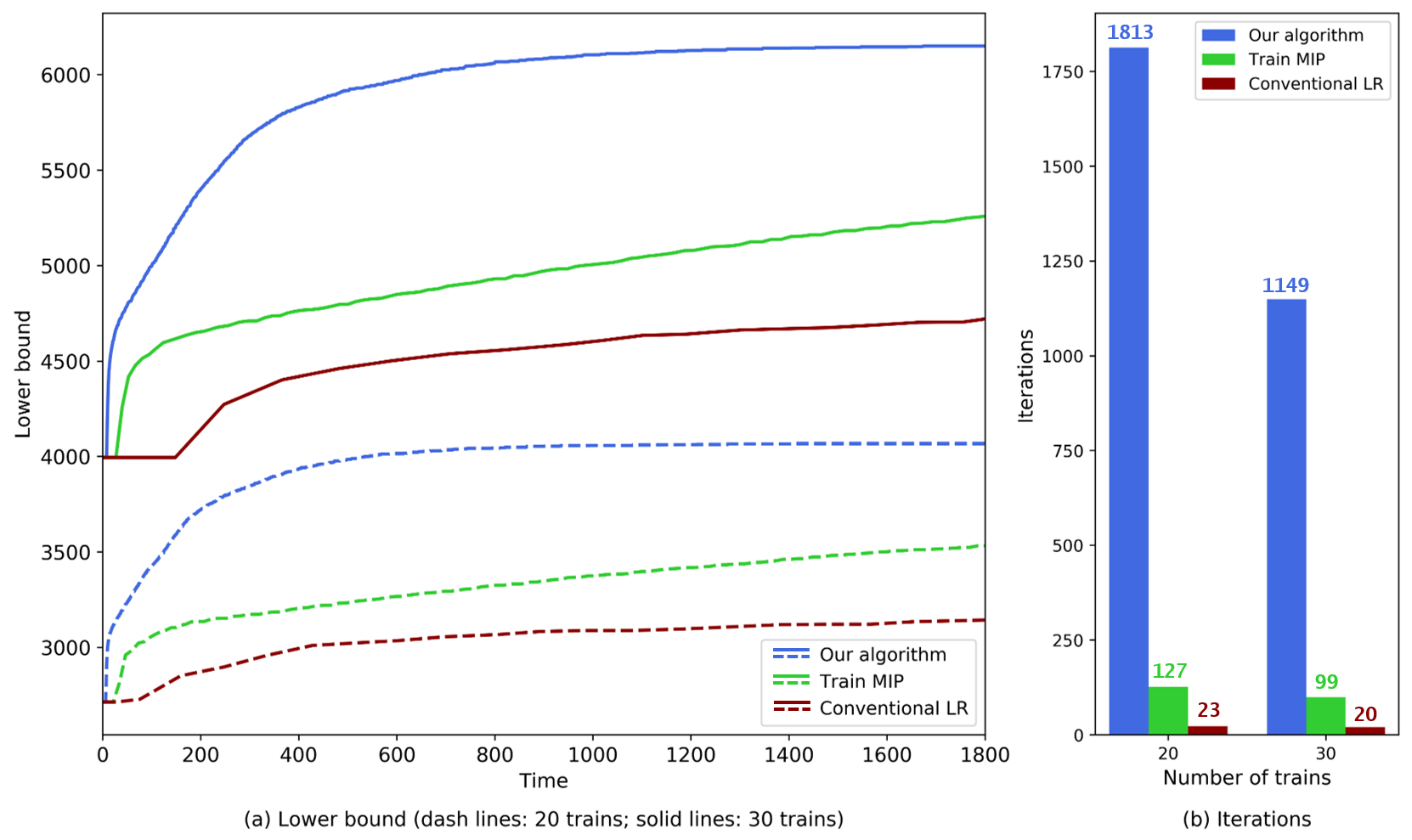}
    \caption{Performance comparison for three methods to solve the relaxed problem.}
    \label{fig:decomposability}
\end{figure}

To validate the effectiveness of the \gls{2-L LR} method, we compare it with two other methods for obtaining lower bounds, namely a {\em conventional single level \gls{LR}} approach and something we term {\em the train MIP approach}. For each approach, we relax the microscopic level capacity constraints~\eqref{headway_route} and~\eqref{siding_capacity} of the full \gls{TLSTN} model into the objective function with the appropriate Lagrange multipliers. The conventional \gls{LR} approach solves the relaxed problem directly. In contrast, the train MIP approach decomposes the relaxed problem into a set of train-independent mixed integer programming models since the objective and constraints in the relaxed problem naturally decompose into train blocks. The mixed integer programming models associated with the train blocks are solved independently. As shown in Fig.~\ref{fig:model_structure}-b), the proposed \gls{2-L LR} approach decomposes the relaxed problem into a master problem and a subproblem. The master problem is further decomposed into a set of train-specific blocks that consider an aggregated penalty for each macroscopic arc. Fig.~\ref{fig:decomposability} compares the lower bounds obtained by the three different methods using a time discretization of five seconds. In each case, the train MIP approach and the proposed \gls{2-L LR} approach can solve the relaxed problem faster in each iteration compared to the conventional \gls{LR} algorithm, resulting in improved convergence within a given time limit. The proposed \gls{2-L LR} is clearly superior to the alternative methods. This is primarily because we solve the master problem only at the macroscopic level, so the solving time required for the relaxed problem in each iteration is significantly reduced. One can observe considerable differences between the total number of iterations for three methods in Fig.~\ref{fig:decomposability}-b).

\begin{figure}[h!]
    \centering
    \includegraphics[scale=0.43]{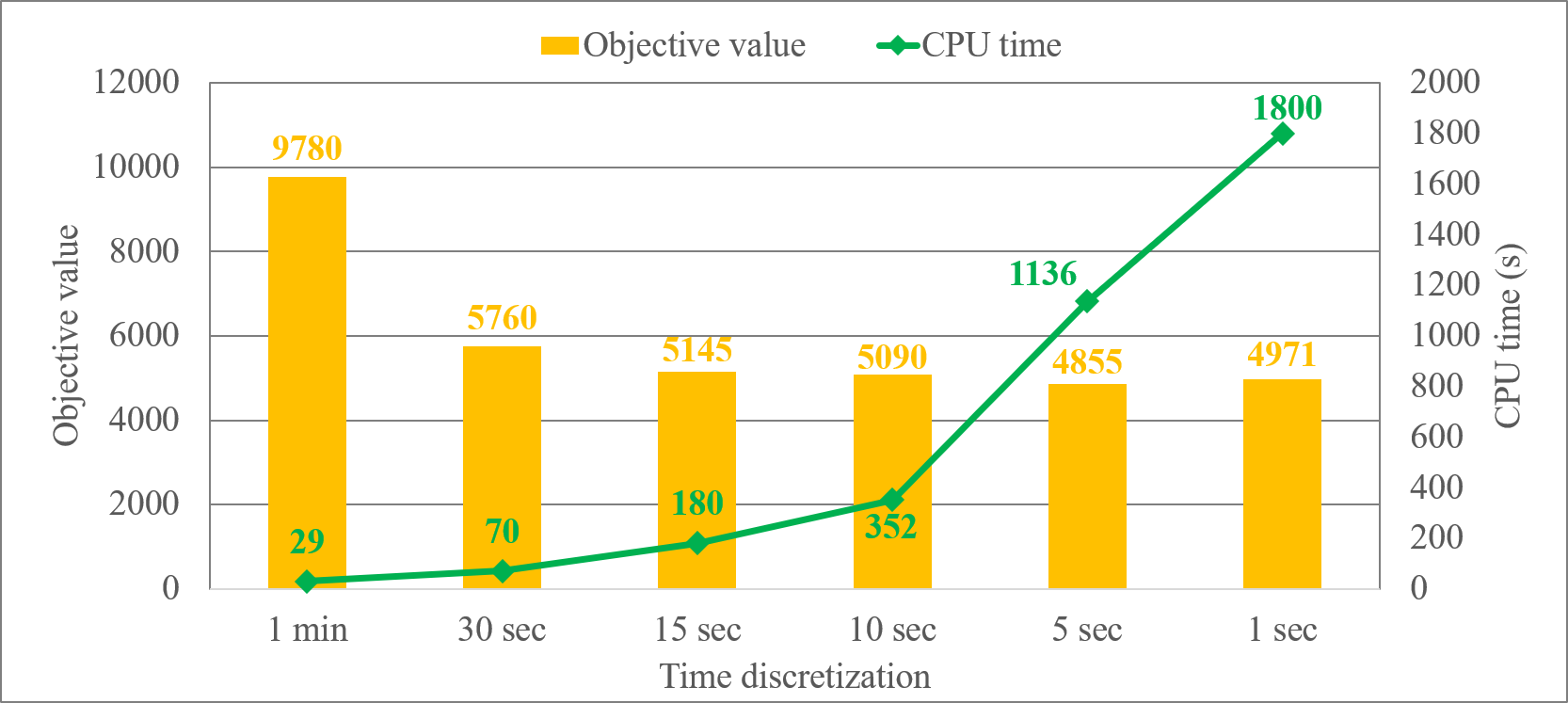}
    \caption{Performance comparison for the time discretization.}
    \label{fig:time_granularity}
\end{figure}

To investigate the effects of time precision on the performance of the proposed algorithm, Fig.~\ref{fig:time_granularity} illustrates the CPU running times and objective values for different values of the time discretization given an instance with 20 trains. The proposed algorithm terminates when the number of iterations exceeds 500 or the solving time exceeds 1800 seconds. When the time discretization is five seconds, the objective value is 4855, which is an improvement of 50\qin{.36}\% compared to using a time discretization of one minute. However, the solution time for the former case is much longer, needing 1136 seconds. One can see that the solution quality improves with finer time discretization, but the solution time also increases. One can also see that the instance with a time discretization of one second fails to yield an improved solution, despite its substantially long solution time. A finer time discretization leads to a larger model with a greater number of arcs and vertices in the constructed space-time network. As a result, the time needed per iteration increases sharply, leading to fewer iterations within the given time limit. Therefore, there is a trade-off between the benefits of solution quality and the solution time required when selecting the time precision for the space-time network.

\subsection{Tests on a large-scale station}

\begin{figure}[h!]
    \centering
    \includegraphics[scale=0.56]{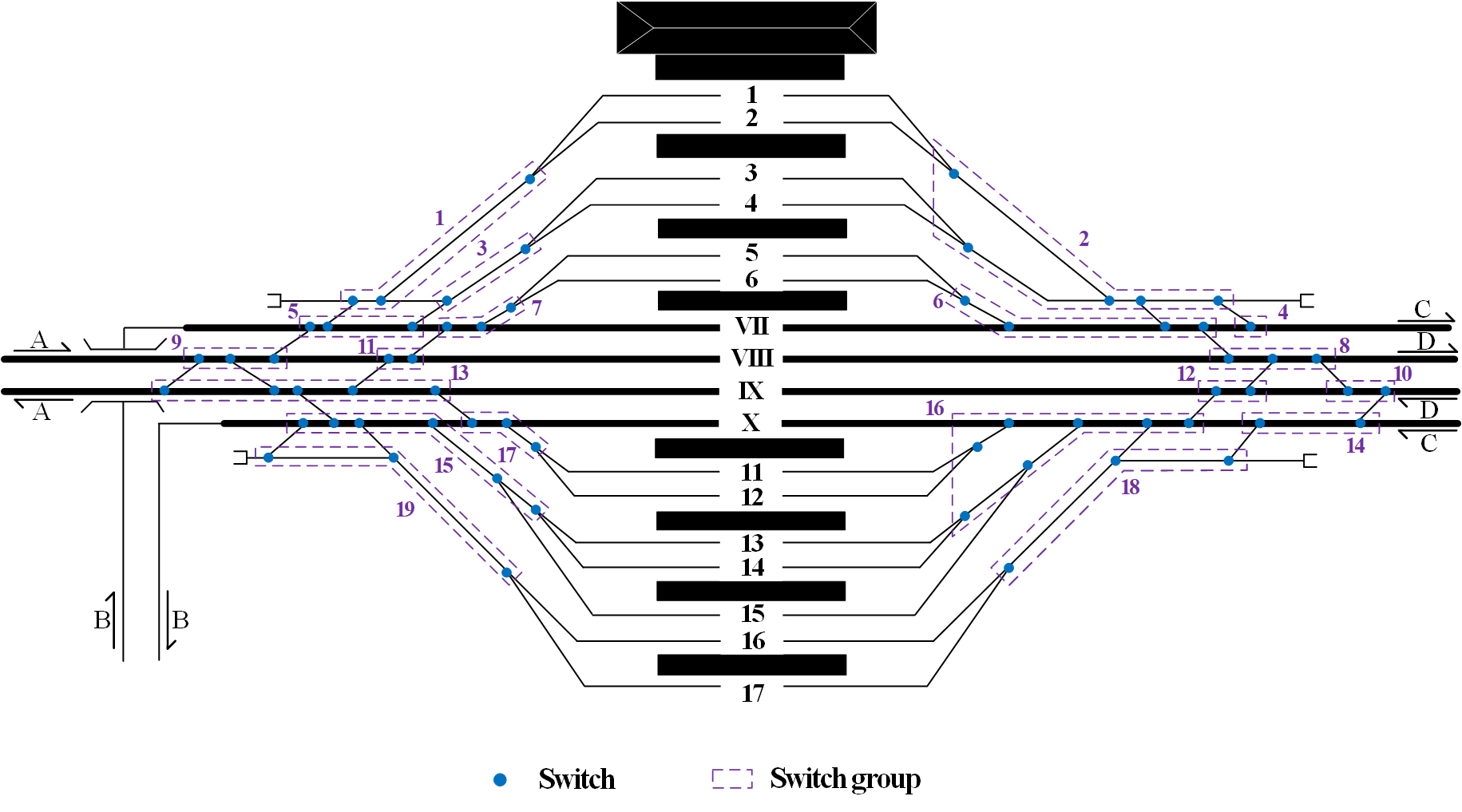}
    \caption{Station layout in large case.}
    \label{fig:station_large_case}
\end{figure}

The large-scale station shown in Fig.~\ref{fig:station_large_case} connects two railway lines and consists of 13 siding tracks and four mainline tracks. This station is located on the busiest Beijing-Shanghai high-speed railway line, and its station layout is commonly used in China. Track~1 at this station is typically occupied by dedicated passenger electric multiple units (EMUs) for the entire day, so this track is not considered in our analysis. We construct the desired timetable based on 287 trains at this station between 4:30 and 00:30 (the next day) for a typical weekday. The headway required for using the same resource is set to 30 seconds. 

\subsubsection{Results analysis at the tactical level}

\begin{figure}[h!]
    \centering
    \includegraphics[scale=0.68]{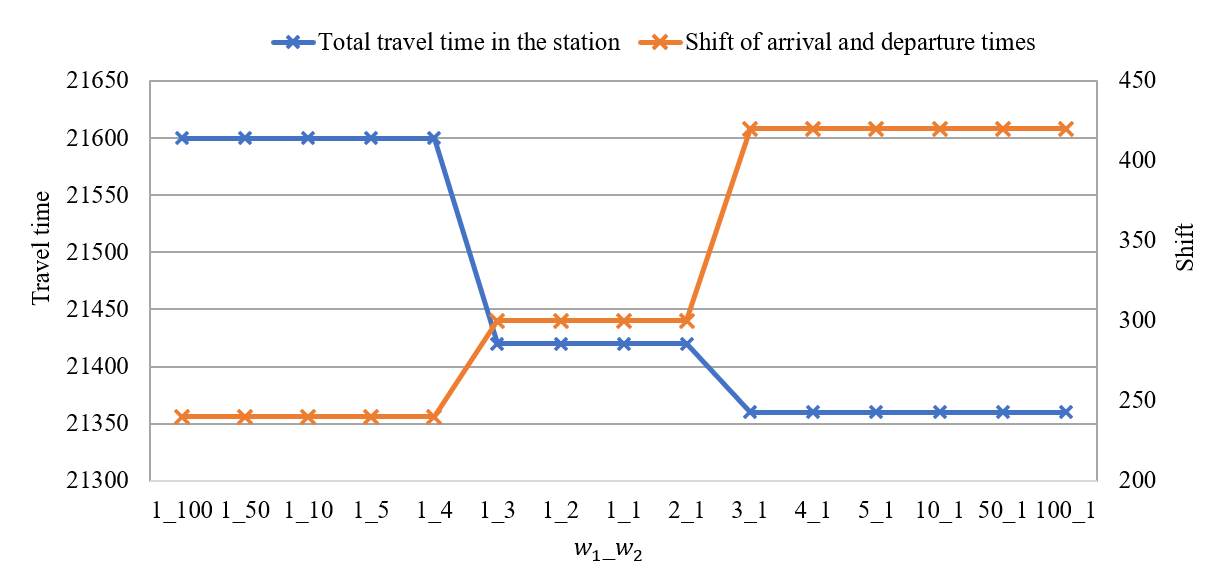}
    \caption{\qin{Trade-off between travel time and shift time.}}
    \label{fig:trade_off}
\end{figure}

\qin{We \rl{first} investigate how the travel time and shift time vary as the values of the corresponding weights in the objective function change. Fifteen 
cases with different combinations of $w_1$ and $w_2$ are solved using the proposed algorithm. We consider 50 trains in every experiment. The maximum solving time and maximum number of iterations are set to 1800 seconds and 500, respectively. Comparisons are shown in Fig.~\ref{fig:trade_off}. The numbers on the horizontal axis correspond to the values of the two weights. For example, 2\_1 indicates that the values of $w_1$ and $w_2$ are set to 2 and 1, respectively. First, we fix the value of $w_1$ to 1. As the value of $w_2$ decreases from 100 to 4, the total shift time remains constant at 240. As $w_2$ decreases further, the total shift time increases to 300 and remains the same until the value of $w_2$ reaches 1. We then fix $w_2$ and increase $w_1$ to 2, resulting in the same value of total shift time as the previous experiment. However, when $w_1$ is set to 3, the shift time increases to 420 and remains the same as $w_1$ increases to 100. In general, the total shift time has a negative correlation with the ratio of $w_1$ to $w_2$, whereas the total travel time is positively correlated with this ratio. Planners can determine specific sets of values according to their preferences. In our large-sized cases, we set the values of $w_1$ and $w_2$ to 1.}

To test the performance of the \gls{2-L LR} algorithm for solving the train platforming problem in practice, we apply it to a case with 287 trains for the entire day at a large-scale station from a planning perspective. 
The allowable \qin{shift}s in arrival and departure times are train-dependent and depend on the desired dwell of a train. Overall, we allow a maximum shift in arrival time of ten minutes (either earlier or later than the desired arrival time) and a maximum extra dwell time of five minutes. 
The time granularity is set to 15 seconds. The algorithm terminates when the number of iterations exceeds 100. The final percentage gap between the best upper and lower bounds is only 1.42\%. 

As shown in Table~\ref{performance_comparison_large_case}, we further compare the performance of our proposed algorithm to that of the \gls{LBBD} and Lazy methods on 12 large-scale cases. Additionally, we report the results obtained by the \gls{2-L LR} algorithm with an upper bound generated only in the final iteration while incorporating a random mechanism. 
We consider four train sets. The first set contains 287 trains for the entire day (W287). The second set contains the first 190 trains of the day with many originating trains (E190). The third set contains 190 trains in the middle of the day with fewer originating and terminating trains (M190). The last set contains the last 190 trains with many terminating trains (L190). For each train set, the instance is further divided into three different scenarios, each of which has a different time discretization. Cases 1,4,7,10 have a fine time discretization of 15 seconds, and Cases 2,5,8,11 have a time discretization of 30 seconds. In the remaining cases, the time discretization is one minute. All experiments are conducted with a time limit of 30 minutes \rl{and not a single case terminates within this time limit.} 
In the table, the first column lists the case ID, the train set and the time discretization. As an example, the entry 1-W287-15 corresponds to Case 1, train set W287, and a 15-second time discretization. Columns two to five present the upper and lower bounds, as well as their percentage gaps and the total iteration within the solving time limit, for the proposed method. The next four columns display the results for the proposed method with the upper bound generated in the final iteration.  The other two methods are unable to generate feasible upper bounds for the problem. Therefore, we only present the lower bounds provided by these methods and the differences between their lower bounds and the best lower bounds obtained with the proposed algorithm.


\renewcommand{\arraystretch}{1}

\begin{table}[h!]
\centering
\caption{Performance comparison for the large-sized cases.}
\label{performance_comparison_large_case}

\small{\small{\small{\small{
\begin{threeparttable}
\begin{tabular}
{p{48pt}p{24pt}p{24pt}p{20pt}p{14pt}p{24pt}p{24pt}p{20pt}p{16pt}p{24pt}p{23pt}p{24pt}p{23pt}}
\toprule
& \multicolumn{4}{c}{\gls{2-L LR}}& \multicolumn{4}{c}{\gls{2-L LR}-UB final generate}&\multicolumn{2}{c}{\gls{LBBD}}&\multicolumn{2}{c}{Lazy-linking} \\
\cmidrule(lr){2-5} \cmidrule(lr){6-9} \cmidrule(lr){10-11} \cmidrule(lr){12-13}
Case&UB&LB&Gap&Iter&UB1&LB1&Gap1&Iter1&LB2&$\Delta{\text{LB}2}$&LB3&$\Delta{\text{LB}3}$\\
\midrule
1-W287-15&201720&197951&1.87\%&17&201510&198164&1.67\%&23&197175&0.50\%&197655&0.26\%\\
2-W287-30&205080&201600&1.70\%&66&204930&201805&1.52\%&107&200700&0.55\%&200970&0.41\%\\
3-W287-60&236660&213395&9.83\%&241&228280&213941&6.28\%&416&212340&0.75\%&212400&0.72\%\\
\midrule
4-E190-15&121485&119992&1.23\%&25&121530&120095&1.18\%&38&119460&0.53\%&119715&0.32\%\\
5-E190-30&123450&122102&1.09\%&114&123420&122245&0.95\%&193&121440&0.66\%&121590&0.54\%\\
6-E190-60&133140&129558&2.69\%&459&132840&129697&2.37\%&608&128700&0.77\%&128875&0.63\%\\
\midrule
7-M190-15&93090&91791&1.40\%&39&93090&91840&1.34\%&52&91140&0.76\%&91230&0.66\%\\
8-M190-30&95100&93899&1.26\%&162&94920&94019&0.95\%&229&93180&0.89\%&93240&0.83\%\\
9-M190-60&104760&101541&3.07\%&485&104100&101764&2.24\%&724&100500&1.24\%&100630&1.11\%\\
\midrule
10-L190-15&125235&122763&1.97\%&26&125235&122826&1.92\%&35&121920&0.74\%&122190&0.52\%\\
11-L190-30&127800&125302&1.95\%&105&127890&125433&1.92\%&145&124470&0.77\%&124650&0.62\%\\
12-L190-60&148180&133712&9.76\%&428&139200&133993&3.74\%&580&132360&1.22\%&132484&1.13\%\\
\bottomrule
\end{tabular}
\begin{tablenotes}
\footnotesize
\item $\text{Gap}=\frac{\text{UB}-\text{LB}}{\text{UB}}*100\%$; $\text{Gap1}=\frac{\text{UB1}-\text{LB1}}{\text{UB1}}*100\%$; $\Delta{\text{LB2}}=\frac{\text{LB2}-\text{LB1}}{\text{LB1}}*100\%$; $\Delta{\text{LB3}}=\frac{\text{LB3}-\text{LB1}}{\text{LB1}}*100\%$; $^*$terminated due to computer memory issues
 \end{tablenotes}
 \end{threeparttable}}}}}
\end{table}

\renewcommand{\arraystretch}{1}

The proposed method is effective at obtaining upper bounds within the time limit for the four cases with a time discretization of 15 seconds. The two upper bound generation methods in our proposed method lead to differences in the lower bounds. This is because generating an upper bound only in the final iteration allows more time, leading to more iterations, to update the lower bounds within the solving time limit. When compared to the results for the small cases, the upper bounds generated in the final iteration are generally better than those generated iteratively. This is because \rl{for} small cases with many conflicts, a random strategy may lead to poor results. \rl{It is also worth mentioning that a larger number of platforms} \qin{and \rl{fewer} trains per unit time leads to \rl{fewer} potential conflicts  \rl{for the} large-sized cases. Therefore, the \rl{optimality} gaps of \rl{the} large-sized cases are smaller than \rl{those of the} small-sized cases 12-17.} The best upper bounds for cases with finer time discretizations are smaller than the tightest lower bounds of corresponding cases with coarser time discretizations, which indicates an improvement in the optimal objective value as the time accuracy decreases. This can be attributed to the longer occupation time of track resources by trains with coarser time discretizations, resulting in more train cancellations in the upper bound. Regarding the train sets, the upper bounds for the cases with the M190 train set are smaller than the lower bounds for E190 and L190. This is because fewer trains originate from or terminate at the station during the planning horizon in the former cases, resulting in less total required dwell time. Furthermore, it should be noted that the \gls{LBBD}-based and Lazy methods are unable to provide upper bounds for large-scale cases. Additionally, these two methods provide slightly weaker lower bounds. 

Across all instances, the proportion of Lagrange multipliers greater than zero accounts for approximately 10\% of the total count. Therefore, the incorporation of the dynamic resource pool significantly reduces the number of Lagrange multipliers that must be considered, ultimately enhancing the performance of the approach. To demonstrate that the model does not have a large number of unnecessary routes, we conducted tests in which trains were only allowed 80\% of their current time shifts. The results show that as the allowable shift decreases (resulting in a decrease in the total number of feasible train paths), the objective function worsened by 0.02-6.65\%. This indicates that there is a benefit in using the specified time shifts.



\subsubsection{Results analysis at the operational level}

Various perturbations can easily influence traffic operations in a high-speed railway system. For example, trains may not arrive or depart from a station as scheduled or some resources may become unavailable. In such cases, the train dispatcher must adjust the track allocation at the station within a short time limit according to the availability of resources and the practical train operation considerations.

Given a feasible platform plan from Table~\ref{performance_comparison_large_case}, we construct two scenarios to highlight the performance of our proposed model and algorithm for solving the train platforming problem at the operational level. First, delays to trains arriving at the station between 15:30 and 17:30 are randomly generated within [0 minute, 10 minutes]. The maximum allowable arrival shift for \rl{each} delayed train is set to \rl{five minutes more than the train's delay}. The allowable arrival shift time for the remaining trains is in the range of [-120 seconds, 120 seconds]. Trains are allowed to dwell for at most an additional 120 seconds. In the second scenario, all trains arrive as scheduled. However, Track 11 is unavailable from 15:30. 
The maximum solution time allowed for our proposed method in these two scenarios is 30 seconds, which is acceptable for operational requirements. 
The proposed algorithm can achieve good results with optimality gaps of only 1.12\% and 1.02\% within the time limit. \qin{\rl{We  refer readers} to https://github.com/zhangqin635/Platform-plan for details.}

\subsubsection{Results with the option for balanced utilization}

In this section, we discuss the results with the option for balancing the use of siding tracks based on constraints~\eqref{balance_constraint}. 
We consider two train sets with 49 trains (between 9:00 and 12:00) and 128 trains (after 15:30) respectively. For each train set, we construct four scenarios with different tolerances of imbalance. We vary the number of paths available for each scenario. The results are 
detailed in~\ref{sec:appendix_balance}. To summarize, decreasing the tolerance for imbalance achieves a more balanced track use. However, the number of cancellations in Cases 1, 4, and 7, where uneven use is acceptable, are fewer than in Case 10, where imbalances are unacceptable. A similar observation can be made for other cases. This is because the excessive pursuit of balanced siding track use may reduce the flexibility of using track resources for trains. As a result, there may be an increase in both canceled train services and objective function values.  There are almost no differences between the standard deviations of instances with the same imbalance tolerance. However, fewer train services are canceled when we allow a wider range of \qin{shift}s in train arrival time or additional dwell time because trains have more routing options.  

\section{Conclusions and future research}
\label{sec:conslusion}
In this paper, we propose a mathematical model for the train platforming problem at a busy high-speed railway station based on a specific two-level space-time network with differences in the modelling precision of the infrastructure. Additionally, we consider different configurations of the interlocking system. In the proposed two-level model, basic operational constraints are enforced in the macroscopic level and complex capacity requirements are guaranteed by constraining resource utilization in the microscopic level. We design a \gls{2-L LR} algorithm to dualize capacity constraints into the objective function and derive a master problem and a subproblem for the relaxed problem. The main difference between our proposed approach and the \gls{LBBD} algorithm is that we evaluate the feasibility of the master solution and penalize the overuse of resources, instead of generating and adding cuts that eliminate macroscopic level solutions. Moreover, we transform Lagrange multipliers in an aggregated manner and decompose the master problem into a set of train-specific blocks, which are independently solved using 
dynamic programming. 
This novel approach can be applied to general two-level problems where one level controls resource utilization and the other enforces utilization feasibility.

In the computational experiments, the proposed formulation and the solution framework are applied to several instances based on a small virtual station and a hub station on one of the busiest Chinese high-speed railway lines to verify their effectiveness. For the small cases, we present quantitative comparisons of two common interlocking configurations described in the literature, namely route-locking, route-release and route-locking, sectional-release. The latter can fully utilize station capacity. Additionally, we compare the performance of our proposed algorithm to that of four other methods, namely the \gls{LBBD} method with no-good cuts, a Lazy approach, a direct CPLEX solve, and an iterative implementation of \gls{LBBD} with no-good cuts. Our \gls{2-L LR} algorithm can obtain superior results with significantly shorter computation time. For the large daily cases on a tactical planning level, we obtain a near-optimal solution with average optimality gaps of around 2\%. We compare the performance of our approach to that of the \gls{LBBD} and Lazy methods on 12 large-scale cases. The results reveal that the \gls{2-L LR} algorithm can provide better upper and tighter lower bounds. Additionally, we test our model and solution method in two disruption scenarios at the operational level. Solutions with optimality gaps of 1.12\% and 1.02\% can be obtained in 30 seconds. We also test the consideration of even track use and the results show that requiring a very balanced platform schedule may substantially limit the routings that are possible and thereby increase the number of train cancellations.

Various applications of two-level network structures in variants of train scheduling and rescheduling problems can be explored. For instance, we can capture different features \rl{at the} different levels when jointly optimizing train timetabling and platforming. Another interesting direction for future research is to further improve the algorithm. 
Various strategies, such as efficient heuristics and effective network reductions, can be explored to generate tighter upper and lower bounds.

\section*{Declarations of interest}
Declarations of interest: none

\section*{Acknowledgement}

This research was supported by the Joint Funds of the National Natural Science Foundation of China (No.U2034208).

\bibliographystyle{elsarticle-harv}
\bibliography{EJOR-temp-supplementary2}

\appendix

\section{Description of the dynamic multiplier generation} 
\label{sec:appendix_dp}

A detailed description of the dynamic multiplier generation is given in~Algorithm~\ref{alg:dynamic_multiplier}. 

\begin{algorithm}[h!]
\caption{Dynamic multiplier generation}\label{alg:dynamic_multiplier}
\begin{algorithmic}
\State \textbf{Input:} $\left\{X_{\text{LB}}^{f}\right\} \forall f\in F$, current iteration $m$, arc pool associated with routes $P_{\text{route}}^m$, microscopic space-time resource pool $M^m$ associated with multipliers, and the value of multipliers for resources in $M^{m}$ 
\State \textbf{Output:}  $P_{\text{route}}^{m+1}$, $M^{m+1}$ and the value of multipliers for resources in $M^{m+1}$ 
\State{\textbf{//Step 1: Initialization}}
\State{Initialize  $M^{m+1} \gets \O$}.
\State{\textbf{//Step 2: Update occupation of trains}}
\State{Initialize occupation of \gls{SG} resources $O_1\gets \O$ and occupation of platform resources $O_2\gets \O$}.
\For{$f\in F$}
  \For{$g\in X_{\text{LB}}^{f}$}
  \State{Update $O_1$ and $O_2$}.
  \EndFor
  \EndFor
\State{\textbf{//Step 3: Update $P_{\text{route}}^{m+1}$  in iteration $m+1$}}
\State{$\left\{P_{\text{route}}^{m+1}\gets P_{\text{route}}^{m}\right\}$}.
\For {$l\in O_1$}
\If {capacity constraints related to $l$ are violated}
\For {$g\in A^a\cup A^d$ associated with $l$}
\State{$\left\{P_{\text{route}}^{m+1}\gets g\right\}$}.
\EndFor
\EndIf
\EndFor
\State{\textbf{//Step 4: Update multipliers should be considered in iteration $m+1$}}
\For {$l\in (O_1\cap M^m)\cup(O_2\cap M^m)$}
\State{Update $\lambda_l^{m+1}=\lambda^{m}_l+\alpha(m)\cdot (\text{occupation}-1)$};
\State{$M^{m+1}\gets l$}.
\EndFor
\For {$l\in (O_1 \cup O_2)\setminus M^m$}
\If {any violation of capacity constraints related to $l$}
\State{Update $\lambda_l^{m+1}=\alpha(m)\cdot (\text{occupation}-1)$};
\State{$M^{m+1}\gets l$}.
\EndIf
\EndFor
\For {$l\in M^m\setminus  (O_1\cup O_2)$}
\State{Update $\lambda_l^{m+1}=\max \left\{0,\lambda_l^{m}+\alpha(m)\cdot (\text{occupation}-1)\right\}$};
\If {$\lambda_l^{m+1}\neq 0$}
\State{$M^{m+1}\gets l$.}
\EndIf
\EndFor
\end{algorithmic}
\end{algorithm}

\section{Consideration of the even track use}
\label{sec:appendix_balance}
The modifications of the solution framework that should be made are listed as follows.

(1) The objective function in the \gls{LR} model

\begin{align}
    \min \qin{Z'_{LR}}&=\qin{Z_{LR}}+\sum_{i \in N^{st}_s}{\lambda_i\left(\sum_{f:g\in A_f}{\sum_{g\in A_f^{a}|d_g=i}{x_g^f}}- \bar{n}_{track}-\bar{\theta}_{track}\right)} \nonumber\\
    &=\qin{Z_{LR}}+\sum_{i \in N^{st}_s}{\lambda_i\sum_{f:g\in A_f}{\sum_{g\in A_f^{a}|d_g=i}{x_g^f}}- \sum_{i \in N^{st}_s}{\lambda_i  \left(\bar{n}_{track}+\bar{\theta}_{track}\right)}}
    \label{z_lr_balance}
\end{align}

(2) The generalized \gls{LR} cost of arcs

\begin{equation}
    \hat{c}_g^f=
    \begin{cases}
    w_1 T_g+\widetilde{c}_g+w_2|\hat{d}_g-T^a_f|+\sum_{l\in B^{st}_f|o_l=d_g,\hat{o}_l\in [\hat{o}_g,\hat{d}_g)}{\lambda_l}+\lambda_i& \forall g \in  A_f^{a},i=d_g,\\
    w_1 T_g+\widetilde{c}_g+w_2|\hat{o}_g-T^d_f|+\sum_{l\in B^{st}_f|o_l=o_g,\hat{o}_l\in [\hat{o}_g,\hat{o}_g+\Delta_s)}{\lambda_l}& \forall g \in  A_f^{d},\\
    w_1 T_g+\qin{\sum_{l\in \phi_g^{ss}}{\lambda_l}}& \forall g \in A_f^{st},
    \end{cases}
    \qquad \forall f \in F\label{cost_arcs_balance}
\end{equation}

(3) The update of multipliers
\begin{equation}
    \lambda_i^{m+1}=\max \left\{0,\lambda_i^{m}+\alpha(m)\cdot\left(\sum_{f:g\in A_f}{\sum_{g\in A_f^{a}|d_g=i}{x_g^f}}- \bar{n}_{track}-\bar{\theta}_{track}\right)\right\}\quad \forall i\in N^{st}_s \label{balance_multiplier_update}
\end{equation}

(4) The heuristic method to obtain the upper bound

When we incorporate the constraints~\eqref{balance_constraint} to balance the use of platform tracks, the heuristic method should be revised slightly to obtain a feasible solution. The trains find the shortest path one by one in the order described in the previous section. In this process, we arbitrarily decide that the siding track will no longer be occupied by any subsequent trains when the number of occupations for this track exceeds the value of $\bar{n}_{track}+\bar{\theta}_{track}$.

The results of the cases that consider different imbalance tolerances are detailed in Table~\ref{result_balance}. The first 
to third columns state the case ID, the number of trains, and values of $\bar{\theta}_{track}$. The following two columns state the additional dwell times and allowed arrival shifts for stopping trains at the station. These two parameters define the number of paths available for each train. The sixth column indicates the standard deviation (STDVA) of the occupation of siding tracks, and the last two column give the number of train cancellations and the respective objective values for each instance.

\begin{table}[h!]
\centering
\caption{Results for 49 and 128 trains considering the balance of utilization.}
\label{result_balance}
\small{\small{\small{
\begin{tabular}
{lllrrrrr}
\toprule
Case ID&\# trains&$\bar{\theta}_{track}$&Extra dwell&Arrival shift&STDVA&Cancellation&Obj.\\
\midrule
1&49&$+\infty$&120&60&5.66&1&30193\\
2&49&$+\infty$&120&120&5.72&0&20958\\
3&49&$+\infty$&300&300&5.83&0&20958\\

4&49&4&120&60&3.6&1&30433\\
5&49&4&120&120&3.67&0&21213\\
6&49&4&300&300&3.66&0&21213\\

7&49&2&120&60&2.67&1&30621\\
8&49&2&120&120&2.52&0&21415\\
9&49&2&300&300&2.52&0&21415\\

10&49&0&120&60&0.78&2&40701\\
11&49&0&120&120&0.87&1&32024\\
12&49&0&300&300&0.39&0&23899\\
\midrule
13&128&$+\infty$&120&300&10.41&0&96795\\
14&128&$+\infty$&300&300&9.71&0&96703\\
15&128&$+\infty$&300&600&10.72&0&96522\\

16&128&4&120&300&4.68&1&105818\\
17&128&4&300&300&4.89&0&99514\\
18&128&4&300&600&4.83&0&99325\\

19&128&2&120&300&3.85&2&115605\\
20&128&2&300&300&3.34&2&114462\\
21&128&2&300&600&3.75&1&107123\\

22&128&0&120&300&1.75&5&140493\\
23&128&0&300&300&1.54&4&132203\\
24&128&0&300&600&1.44&3&124787\\
\bottomrule
\end{tabular}}}}
\end{table}

\section{Supplementary information for small cases}
\label{sec:small case_detail}
We list the number of trains, the number of space-time paths, the number of conflicts between train-path-pair and the numbers of variables and constraints for 17 small cases in Table~\ref{small_case_details}. 

\begin{table}[h!]
\centering
\caption{Supplementary information for small cases in Table~\ref{performance_comparison_small_case}.}
\label{small_case_details}
\small{
\small{\small{\small{
\begin{tabular}
{p{35pt}p{40pt}p{50pt}p{50pt}p{50pt}p{60pt}}
\toprule
Case ID& \# trains&\# paths&\# conflicts&\# variables&\# constraints\\
\midrule
1&4&6800&1771890&5340&3884\\
2&4&6800&1137686&5369&3644\\
3&5&7905&1855852&6491&4629\\
4&5&8500&1628547&6699&4573\\
5&6&9010&1873003&7628&5382\\
6&6&10200&2264515&8040&5637\\
7&7&10115&1889900&8777&6082\\
8&7&11305&2339133&9191&6376\\
9&8&11815&1977557&10094&6924\\
10&8&12410&2354306&10329&7095\\
11&9&13515&2363760&11487&7742\\
12&12&18020&3639990&15288&10020\\
13&15&23120&3954493&19733&12540\\
14&18&28220&4269181&24177&15111\\
15&21&33320&5128661&28377&17671\\
16&24&38420&5200003&32544&20119\\
17&27&41722&5315623&35858&22191\\
\bottomrule
\end{tabular}}}}}
\end{table}

\end{document}